\documentclass[11pt]{article} 

\usepackage[utf8]{inputenc} 
\usepackage{float}
\usepackage{geometry} 
\usepackage{amsmath}
\allowdisplaybreaks
\usepackage{blkarray}
\usepackage{multirow}
\usepackage{graphicx} 
\usepackage{setspace}
\usepackage[toc,page]{appendix}
\usepackage{color}
\usepackage{booktabs} 
\usepackage{array} 
\usepackage{paralist} 
\usepackage{amsfonts}
\usepackage{verbatim} 
\usepackage{subfig} 
\usepackage{amsthm}
\usepackage{color}
\usepackage{fancyhdr} 
\usepackage{sectsty}
\allsectionsfont{\sffamily\mdseries\upshape} 
\usepackage[mathscr]{euscript}
\usepackage[nottoc,notlof,notlot]{tocbibind} 
\usepackage[titles,subfigure]{tocloft} 

    \newtheorem{theorem}{Theorem}
    \newtheorem{lemma}{Lemma}
    \newtheorem{proposition}{Proposition}
    \newtheorem{corollary}{Corollary}

\usepackage{arydshln}
\newtheorem{definition}{Definition}

    \newcounter{example}
    \newenvironment{example}[1][]{\refstepcounter{example}\par\medskip\noindent%
       \textbf{Example~\theexample. #1} \rmfamily}{\medskip}
    
\title{The Value of Misinformation and Disinformation}
\author{Yanling Chang$^a$, Matthew F. Keblis$^a$, Ran Li$^b$, Eleftherios Iakovou$^a$, and Chelsea C. White III$^c$\\
$^a$Texas A\&M University, $^b$Sabre Corporation, and
$^c$Georgia Institute of Technology}
\begin{document}
\maketitle
\doublespacing
\begin{abstract}
Information is a critical dimension in warfare. Inaccurate information such as misinformation or disinformation further complicates military operations.  In this paper, we examine the value of misinformation and disinformation to a military leader who through investment in people, programs and technology is able to affect the accuracy of information communicated between other actors. We model the problem as a partially observable stochastic game with three agents, a leader and two followers. We determine the value to the leader of misinformation or disinformation being communicated between two (i) adversarial followers and (ii) allied followers. We demonstrate that \textit{only under certain conditions}, the prevalent intuition that the leader would benefit from less (more) accurate communication between adversarial (allied) followers is valid.  We analyzed why the intuition may fail and show a holistic paradigm taking into account both the reward structures and policies of agents is necessary in order to correctly determine the value of misinformation and disinformation.  Our research identifies efficient targeted investments to affect the accuracy of information communicated between followers to the leader's advantage.	
\end{abstract}

\section{Introduction}
Both state and non-state actors seek to advance their objectives by leveraging the information environment (IE) defined as “the aggregate of individuals, organizations, and systems that collect, process, disseminate, or act on information” (Joint Publication 1-02).  The recent Battle of Mosul (2016-2017) is representative.  As reported in {\it The Wall Street Journal}, the battle involved just a few thousand untrained militants arrayed against a much larger formally trained force (Spencer 2017). Over the course of more than two years prior to the battle, the militants constructed elaborate defenses including an extensive tunnel system. During the nine months long conflict afterwards, the militants used social media to communicate amongst themselves the exact location and direction of opposing forces, which enabled the militants to move mobile defense assets (everything from snipers to exploding vehicles) to specified locations at specified times to achieve their desired objectives (e.g., mortality of opposing forces).  While information sharing via social media aided the efforts of the militants, communication problems on the other hand hindered the opposing forces.   One particular problem arose in providing fire support to units advancing against militant positions.  Units providing fire support needed information about the ``forward line of troops'' to prevent friendly force casualties, however communicating accurate information between units in the battlefield was very challenging for a variety of reasons including a changing electromagnetic spectrum (Mosul Study Group 2017).\\

The above illustrates the centrality of the IE militarily.  Not surprisingly, leading state actors have long seen dominating the IE as crucial to the success of their efforts in warfare (Department of the Army 1996).  Dominance of the IE involves (i) preserving one's own capabilities and those of allies, and (ii) affecting the capabilities of rivals (Ardis and Keene 2018).  Historically, dominating the IE in open war has to a large degree been about destroying the command, control and communication (C3) infrastructure of an adversary.  However, destruction of C3 infrastructure has become more difficult with the rise of the internet, advances in information technology, widespread availability of wireless communication, and the ubiquity of social media (Department of Defense 2016). Moreover, with advanced technology now more widely available and less costly than it was before, states with limited resources and non-state actors, both previously relegated to the sidelines, have emerged as rivals to traditional powers in the IE.  The result of these trends is that the IE is now a much more competitive space than it was previously and establishing dominance of the IE is becoming more difficult.  Nonetheless, many actors continue to attempt to shape the IE to their advantage, where shaping entails making investments in people, programs and technology that preserve/improve friendly mission-essential information and affect rival mission-essential information (JCOIE 2018).\\

In this paper, we determine the value to an actor (e.g., a military force; referred to as the leader in the rest of the paper) of shaping the IE through investment.  Specifically,  we analyze the value to the leader of modulating information flows between other actors (referred to as followers). We examine two types of scenarios: (i) a non-collaborative scenario, where followers are adversaries of the leader (e.g., groups of militants in the Battle of Mosul; the leader is the opposing military force), and (ii) a collaborative scenario, where followers are allies of the leader (e.g., units of the military force opposing the militants in the Battle of Mosul; the leader is the command element of the military force). \\

In each scenario, the leader's performance is determined by the actions and policies of all agents, whereas the situational understanding of the followers and the decisions they make are affected by information communicated from other followers. The information passing between followers is, however, not necessarily accurate. Inaccurate information is either disinformation or misinformation.  Disinformation is false information spread deliberately to deceive, whereas misinformation is false or incorrect information spread unintentionally (without realizing it is untrue). In the following, we refer to either type of information as distorted. To capture the dynamic nature of warfare, we model the problem as a general-sum partially observable stochastic game (see Hespanha and Prandini 2001, Bernstein et al. 2002, Rabinovich et al. 2003, Hansen et al. 2004, Kumar and Zilberstein 2009) with three agents, a leader and two followers. The research is intended to provide military commanders, defense planners and military analysts with insights pertaining to the value of distorted information in warfare. We further identify efficient targeted investments for the leader in order to shape the IE to the leader's advantage.  \\

\textbf{Related Literature}. The study of distorted information in communication has attracted the attention of both the social and computer science communities.  Within the former, Lewandowsky et al. (2012) reviewed the ways in which information is communicated and distorted information is disseminated within society while Lewandowsky et al. (2013) used two case studies to illustrate how distorted information could be employed to either escalate or de-escalate violent conflicts. Research in the computer science community has focused on modeling the spread of distorted information and developing interdiction strategies to limit the spread of distorted information in social networks. Nguyen et al. (2012) utilized a zero-sum game to study how a defending party can minimize the impact of distorted information spread by an adversary. Budak et al. (2011) evaluated various approaches for limiting the spread of distorted information in a social network, while Zhang et al. (2016) developed several algorithms to detect, with limited budget, distorted information nodes in such networks. A survey of detecting distorted information in social media is provided by Zubiaga et al. (2017).\\

\begin{figure}[h]
	\centering
	\begin{minipage}{0.24\textwidth}
		\centering
		\includegraphics[width=0.70\textwidth]{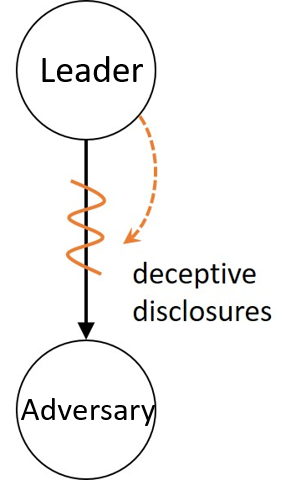}\\
		{(a) In literature and our previous work}
	\end{minipage}\hfil
	\begin{minipage}{0.31\textwidth}
		\centering
		\includegraphics[width=\textwidth]{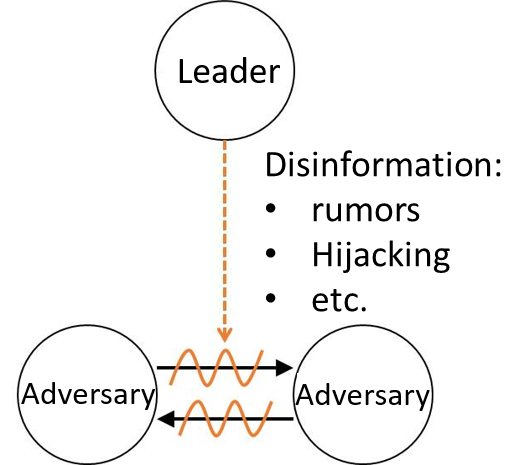}\\
		{(b) Non-cooperative games}
	\end{minipage}\hfil	\begin{minipage}{0.45\textwidth}
		\centering
		\includegraphics[width=\textwidth]{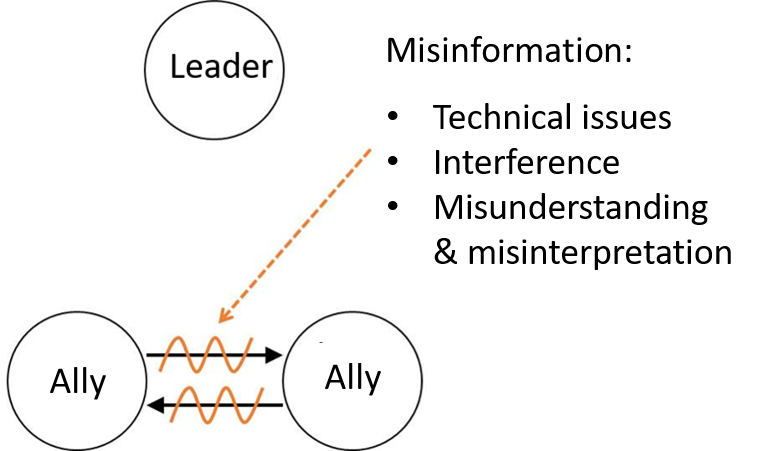}\\
		{(c) Cooperative games}
	\end{minipage}\hfil
	\caption{(a) Distorted information in literature and our previous work vs. distorted information in communication for the non-cooperative (b) and cooperative games (c) studied in this paper}
	\label{difference}
\end{figure}

In the decision analysis and game theory literatures, single-period games involving two agents have been extensively utilized in value of information analysis. For example, Li (2002), Chu and Lee (2006), Meyer et al. (2010), and Leng and Parlar (2009) all examined the positive value of revealing private information, while Kamien et al. (1990) studied problems where rewards could degrade with more accurate information. Lehrer and Rosenberg (2010) further characterized the value-of-information function for zero-sum repeated game with incomplete information. In the security context, Greenberg (1982), Ma et al. (2013) and Chang et al. (2015a, 2015b) all have concentrated on the role of distorted information using \textit{two-agent} games, where the leader manipulates the information (state, actions, etc.) available about himself that is obtainable by the follower (the adversary); see Figure \ref{difference}(a).  A review of this type of information distortion is provided by Merrick et al. (2016).\\

The focus of this paper is quite different from the foci encountered in the existing literature and our previous work Chang et al. (2015a, 2015b) as we study the value of distorted information in \textit{communication} for the situations depicted in Figure \ref{difference}(b,c).  Figure \ref{difference}(b) corresponds to a non-collaborative scenario where an army has the objective of reclaiming an area held by groups of militants.  Each group of militants maintains a position within the area and monitors the activities of the army in the vicinity of its position. As reported in Spencer (2017), each group of militants is able to gather information (number of opposing soldiers, types of weapons, etc.) only about the army units operating in the vicinity of its position.  However, each group of militants can share information with other groups of militants in other positions in order to gain a more complete picture of the army's activities.  Each group of militants uses the information it receives from other militant groups about the army's activities to inform its decisions, e.g., about where to position heavy weapons, move militant soldiers, etc.  In this paper, we determine the value to the army of modulating information communicated between groups of militants (e.g., disinformation through hijacking), in order to guide investment decisions the army can make to affect adversary mission-essential information.\\

Figure \ref{difference}(c) corresponds to a collaborative scenario where a military force, consisting of a command element and combat units, has the objective of retaking an area held by groups of militants.  Each combat unit is cognizant of the progress it is making towards taking a position held by one of the militant groups, however each combat unit has a less than perfect understanding of the progress other combat units are making against other positions.  Each combat unit attempts to communicate to other combat units its position on the battlefield and uses information received from other combat units to inform whether and where to provide fire support to other combat units.  Accurate information is vital given the close proximity of targets to the forward line of troops. However, real-time battlefield communication is challenging (Mosul Study Group 2017).  In urban warfare, misinformation is not uncommon due to technical issues associated with the communication technologies (Fitzgerald et al. 2014) and interference caused by buildings and structures that impede electronic communication (Edwards 2001).  Such communication difficulties combined with the ever-changing nature of the urban environment during warfare have resulted in combat units on the same battlefield having a very different situational understanding, as noted in {\it The Wall Street Journal} (Gordon 2017).  In this paper, we determine the value to the command element of modulating information communicated between combat units, in order to guide investment decisions the command element can make to preserve/improve allied mission-essential information (e.g., through investment in advanced communication technology, development of common understanding of symbols).\\

We next discuss the contributions of our paper.
\begin{enumerate}[(i)]
	\item Our developed modeling framework constitutes a first effort in timely and rigorously attempting to capture and analyze the value of distorted information in communication in modern warfare from the leader's perspective. We model the IE as a general-sum partially observable stochastic game with three agents to capture the dynamic nature of urban warfare. The value of distorted information in both cooperative and non-cooperative scenarios is analyzed.
	\item A holistic paradigm considering both the reward structures and the policy employed by each agent is suggested for analyzing the value of distorted information in communication. Embracing such a holistic view of the reward structure, policy, and accuracy of information is imperative as considering only one of them can lead to erroneous results.  Intuition suggests that the leader's performance should improve (deteriorate) as the quality/accuracy of the information communicated between adversarial (allied) followers declines.  This is widely accepted in the existing literature.  Herein, we demonstrate that such intuition is correct only under certain conditions. We also discuss why the intuitive understanding of distorted information can be false.
	\item Mathematically, we show the leader's value function can be represented as a power series of information accuracy communicated between two followers. Thus, the value function of the leader is infinitely differentiable in communication quality. We identify several sets of conditions under which the intuitive understanding of distorted information is valid.  
	\item In situations where the leader would benefit from modulating the information communicated between followers, we analyze the connection between affecting information accuracy and efficient investment. This analysis provides guidance to a leader with limited resources about how to invest in the modulation of communication between followers.
\end{enumerate}

The rest of the paper is organized as follows. Section \ref{problemstatementsec} formulates the problem as a one-leader, two-follower partially observable stochastic game. We state the modeling assumptions in Section \ref{ModelAssumptions} and the model formulation in Section \ref{modelformulation}. Section \ref{communicationQuality} presents the assessment of the value of distorted information. We first analyze the best response behavior of the followers, and discuss its implication on the value of distorted information. Subsequently, we represent the leader's value function as a power series of communication quality in Section \ref{sectionPowerSeries}. Section \ref{SectionImpact} analyzes the impact of the reward structures and the polices employed by each agent on the value of distorted information. Section \ref{othercases} examines the cases where the analysis in Section \ref{sectionPowerSeries}-\ref{SectionImpact} is not applicable. We analyze why the intuitive understanding of distorted information could be invalid and discuss the implications on investment decisions in Section \ref{Implication}. Section \ref{conclusion} concludes the paper where we summarize our work and discuss future research directions.

\section{The Game Setup}
\label{problemstatementsec}
In this section, we model the problem as a partially observable stochastic game involving three agents: a leader, follower 1, and follower 2. The two followers are adversaries of the leader in a non-cooperative game as depicted in Figure \ref{difference}(b) and allies of the leader in a cooperative game as depicted in Figure \ref{difference}(c).
\subsection{Model Assumptions}
\label{ModelAssumptions}
We assume that (i) each agent has perfect knowledge of the game setup and (ii) before the game, the leader chooses his policy first and then the followers simultaneously select their own policies with complete knowledge of the policy chosen by the leader. In cooperative scenarios, the leader may directly inform allies of his policies and of the game setup, while in non-cooperative scenarios, the adversaries may have spent a substantial amount of time and effort in acquiring knowledge of the game and in understanding the leader's strategies before making any decisions (e.g., see \textit{Joint Publication 3-13}).\\

We remark that knowing an agent's policy is not equivalent to knowing the agent's action at each decision epoch in sequential decision making. The action of an agent is determined by the policy of the agent in concert with the information the agent possesses. The information possessed by an agent, however, is private and is not accessible to other agents. We study how the real-time passing of distorted information between the two followers can affect the decision making of the followers and hence the performance of the leader.

\subsection{Model Formulation}
\label{modelformulation}
\begin{figure}[h]
	\centering
	\includegraphics[width=0.55\textwidth]{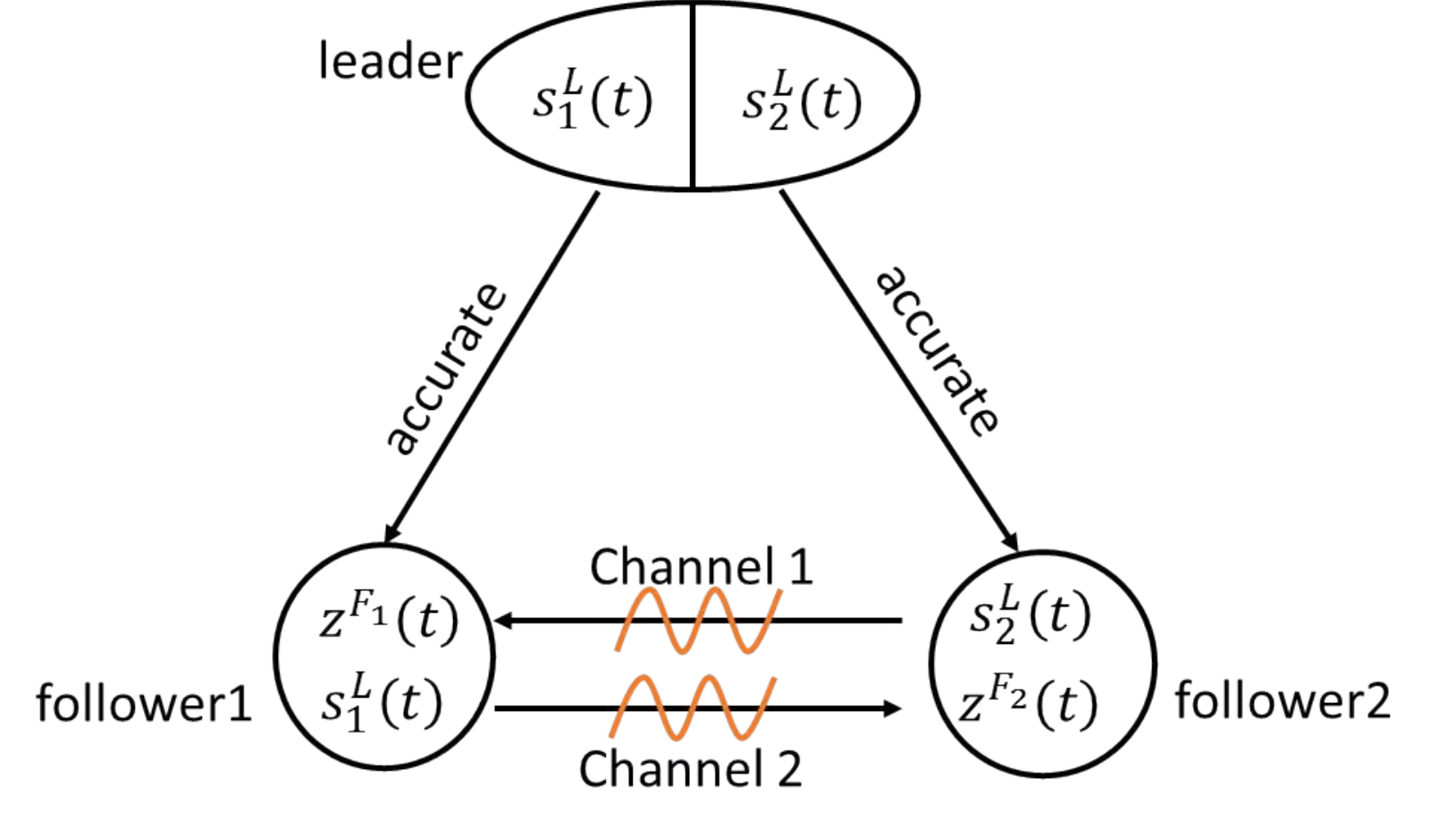}
	\caption{Communication between two followers}
	\label{scheme}
\end{figure}

The decision epochs are $t=0,1,2,...$, and let $\{s^k(t), t=0,1,...\}$, $\{\bar{z}^k(t), t=1,2,...\}$, and $\{a^k(t), t=0,1,...\}$ be the state, observation, and action processes for agent $k\in \{L=\text{leader}, F_1 = \text{follower 1}, F_2 = \text{follower 2} \}$, each having finite state, observation, and action spaces $S^k$, $\bar{Z}^k$, and $A^k$, respectively. Observation $\bar{z}^k(t)$ received by agent $k$ contains possibly distorted information of the other two agents' states.  A follower can independently collect partial information about the leader and share information with the other follower through a communication channel, e.g., a social media application. We assume: (i) the leader's state $s^L(t)$ has two components: $s^L(t) = (s^L_1(t), s^L_2(t))$; the state space of the leader is given by $S^L = S^L_1 \times S^L_2$; (ii) follower $i \in \{1,2\}$ can perfectly observe the leader's partial state $s^L_i(t)$ at time $t$, whereas follower $i`s$ knowledge of the leader's partial state $s^L_{3-i}(t)$ results from (possibly distorted) information $z^{F_i}(t) \in Z^{F_i}$ obtained from follower $3-i$ via communication channel $i$ (see Figure \ref{scheme}).  Assume $Z^{F_i} = S^L_{3-i}, i \in \{1,2\}$. We say the two followers accurately share information with each other at any time to coordinate their actions if $z^{F_i}(t) = s^L_{3-i}(t), \forall t$.  However, the existence of distorted information in communication may result in $z^{F_i}(t) \neq s^L_{3-i}(t)$.

As our focus is on the distorted information in communication between two followers, we assume that: (i) each follower $i \in \{1,2\}$ has complete knowledge of his own state and of the state of the other follower, together denoted as $s^F(t)=\{s^{F_1}(t), s^{F_2}(t)\}$, hence $\bar{z}^{F_i}(t) =\{s^L_i(t),z^{F_i}(t), s^{F_{3-i}}(t)\}$; and (ii) the leader can perfectly observe $s^F(t)$, hence $\bar{z}^L(t) = s^F(t)$.  We remark that all the results of this paper can be directly extended, albeit with more complicated notation, to the case where each follower does not have complete knowledge of the state of the other follower and where the leader's observation of $s^F(t)$ is less than perfect.\\

In the following, let $s(t) = \{s^L(t), s^{F_1}(t), s^{F_2}(t)\}$, $\bar{z}(t) = \{\bar{z}^L(t), \bar{z}^{F_1}(t), \bar{z}^{F_2}(t)\}$ and $a(t) = \{a^L(t), a^{F_1}(t), a^{F_2}(t)\}$.  Furthermore, let $Q=\{P(\bar{z}(t)|s(t))\}$ be the communication matrix (also called the observation matrix in the POMDP literature). Then $Q(\bar{z}(t)|s(t))=Q(z(t)|s(t))=\prod_{i=1}^2 Q^i(z^{F_i}(t)|s^L_{3-i}(t))$, where $Q^i$ is the communication matrix for channel $i, i \in \{1,2\}$ and $z(t)=\{z^{F_1}(t),z^{F_2}(t)\}$. The higher communication quality of a channel, the lower the probability of distorted information being passed in an information exchange. Let $\epsilon_i \geq 0$ represent the probability of distorted information in communication channel $i$. We thus define the communication matrix for channel $i \in \{1,2\}$ as:
$$Q^i(z^{F_i}(t) |s^{L}_{3-i}(t))=
\begin{cases}
1- \epsilon_i &\mbox{ if } z^{F_i}(t) = s^{L}_{3-i}(t) \\
\epsilon_i \cdot \sigma_{s^{L}_{3-i},z^{F_i}} &\mbox{ if } z^{F_i}(t) \neq s^{L}_{3-i}(t),	
\end{cases}$$
where $\sigma_{s^{L}_{3-i},z^{F_i}} \geq 0,  \sigma_{s^{L}_i,s^{L}_i}=0$, $ \sum_{z^{F_i} \in Z^{F_i}}\sigma_{s^{L}_{3-i},z^{F_i}}  = 1$ for all $s^{L}_i, i \in \{1,2\}$. We further define $||Q||=\max_j \sum_{z \in Z}|q_{jz}|$ for a matrix $Q$. Clearly, $||Q^i - I|| = 2\epsilon_i, i \in \{1,2\}$, where $I$ is the identity matrix. The communication quality of channel $i$ is considered higher if its communication matrix is given by $Q^i(\epsilon_i)$ rather than $Q^i(\epsilon_i')$ where $\epsilon_i < \epsilon_i'$.
We remark that the definition of $Q^i$ is equivalent to the definition of the observation matrix in Ortiz et al. (2013) for POMDPs. Another definition of information accuracy in communication will be discussed later in Subsection \ref{othercases}. The conditional probability $p_{ij}(\bar{z},a) = P[\bar{z}(t+1)=\bar{z}, s(t+1)=j|s(t)=i, a(t)=a]$ is assumed given, and let $P(\bar{z},a)$ be the sub-stochastic matrix $\{p_{ij}(\bar{z},a)\}$.\\

Let $\zeta^k(t,\tau) = \{s^k(t),...,s^k(t-\tau+1), \bar{z}^k(t),..., \bar{z}^k(t-\tau+1), a^k(t-1),...,a^k(t-\tau)\}$ be the information history at time $t$ of finite length $\tau$ for agent $k$, and $\bar{\zeta}$ be such that $\zeta^k(t+1,\tau) = \bar{\zeta}(\bar{z}^k(t+1), s^k(t+1), a^k(t), \zeta^k(t,\tau))$. Let $\zeta^{l,p}(t,\tau) = (\zeta^l(t,\tau), \zeta^p(t,\tau)), \zeta(t,\tau) = (\zeta^L(t,\tau), \zeta^{F_1}(t,\tau), \zeta^{F_2}(t,\tau))$, and $y^k(t) = \{P(\zeta^{l,p}(t,\tau)|\zeta^k(t)) \}$ for $l \neq p \neq k$, where $y^k(t)$ is a ``belief'' array indicating what agent $k$ infers about the other two agents' information histories $\zeta^{l,p}(t,\tau), l \neq p \neq k, l, p, k \in \{L, F_1, F_2\}$. Assume the initial belief $y^k(0) = \{P(\zeta^{l,p}(0,\tau))\}$ is given. Let $\zeta^k(t) = \{\bar{z}^k(t),...,\bar{z}^k(1), s^k(t),...,s^k(0),a^k(t-1),...,a^k(0), y^k(0)\}$ when $t \geq 1$, hence, $\zeta^k(t) = \{\bar{z}^k(t),s^k(t),a^k(t-1),\zeta^k(t-1)\}$. A policy $\delta^k$ is called a \textit{perfect-memory} policy if agent $k$ selects $a^k(t)$ on the basis of complete information history $\zeta^k(t)$; it is called a \textit{finite-memory} policy if agent $k$ selects $a^k(t)$ on the basis of finite information history $\zeta^k(t,\tau)$. In this paper, we assume the policies employed by the leader and followers are finite-memory policies.  For simplicity, we define $\delta=(\delta^L, \delta^{F_1}, \delta^{F_2})$, and $\delta(\zeta(t,\tau)) = (\delta^L(\zeta^L(t,\tau)), \delta^{F_1}(\zeta^{F_1}(t,\tau)), \delta^{F_2}(\zeta^{F_2}(t,\tau)))$ = $P(a(t)|\zeta(t,\tau))$, where $P(a(t)|\zeta(t,\tau)) = \prod_{k \in \{L,F_1,F_2\}} P(a^k(t)|\zeta^k(t,\tau))$. We remark that $\zeta^k(t,\tau)$ is \textit{private} information of agent $k$ and it is only known by the agent. Hence, agent $l$ can only \textit{infer} agent $k$'s action even if agent $l$ knows the policy of agent $k$, $k \neq l$. For simplicity, we further denote:
\begin{enumerate}[(i)]
	\item  $\bar{\zeta}(s(t+1),s^{L}_2(t+1),s^{L}_1(t+1),\zeta(t,\tau),\delta(\zeta(t,\tau))) = \{ \bar{\zeta}(s(t+1),\zeta^L(t,\tau),\delta^L(\zeta^L(t,\tau))),\bar{\zeta}(s(t+1),\zeta^{F_1}(t,\tau),\delta^{F_1}(\zeta^{F_1}(t,\tau))),\bar{\zeta}(s(t+1),\zeta^{F_2}(t,\tau),\delta^{F_2}(\zeta^{F_2}(t,\tau)))\}$ (each agent receives accurate information at time $t+1$);
	\item  $\bar{\zeta}(s(t+1),z^{F_1}(t+1),s^{L}_1(t+1),\zeta(t,\tau),\delta(\zeta(t,\tau))) = \{ \bar{\zeta}(s(t+1),\zeta^L(t,\tau),\delta^L(\zeta^L(t,\tau))),\bar{\zeta}(s^{L}_1(t+1), z^{F_1}(t+1),s^{F}(t+1),\zeta^{F_1}(t,\tau),\delta^{F_1}(\zeta^{F_1}(t,\tau))),\bar{\zeta}(s(t+1),\zeta^{F_2}(t,\tau),\delta^{F_2}(\zeta^{F_2}(t,\tau)))\}$ (follower $1$ receives inaccurate information but others' information are accurate at time $t+1$);
	\item  $\bar{\zeta}(s(t+1),s^{L}_2(t+1),z^{F_2}(t+1),\zeta(t,\tau),\delta(\zeta(t,\tau))) = \{ \bar{\zeta}(s(t+1),\zeta^L(t,\tau),\delta^L(\zeta^L(t,\tau))),\bar{\zeta}(s(t+1),\zeta^{F_1}(t,\tau),\delta^{F_1}(\zeta^{F_1}(t,\tau))),\bar{\zeta}(z^{F_2}(t+1),s^{L}_2(t+1), s^{F}(t+1),\zeta^{F_2}(t,\tau),\delta^{F_2}(\zeta^{F_2}(t,\tau)))\}$ (follower $2$ receives inaccurate information but others' information are accurate at time $t+1$);
	\item  $\bar{\zeta}(s(t+1),z^{F_1}(t+1),z^{F_2}(t+1),\zeta(t,\tau),\delta(\zeta(t,\tau))) = \{ \bar{\zeta}(s(t+1),\zeta^L(t,\tau),\delta^L(\zeta^L(t,\tau))),\bar{\zeta}(s^{L}_1(t+1),z^{F_1}(t+1),s^{F}(t+1), \zeta^{F_1}(t,\tau),\delta^{F_1}(\zeta^{F_1}(t,\tau))),\bar{\zeta}(z^{F_2}(t+1),s^{L}_2(t+1),s^{F}(t+1), \zeta^{F_2}(t,\tau),\delta^{F_2}(\zeta^{F_2}(t,\tau)))\}$ (both followers receive inaccurate information at time $t+1$).
\end{enumerate}

Let $r^k(s(t), a(t))$ be the scalar reward received by agent $k$ at epoch $t$ with state $s(t)$ and action $a(t)$. Let $v^k_{\delta}(\zeta^k(0))$ be the agent $k$'s value function under a leader policy $\delta^L$, follower 1's policy $\delta^{F_1}$,  and follower 2's policy $\delta^{F_2}$, with the communication matrix $Q(\epsilon_1, \epsilon_2)$.  As the exact ending period of a war front is unknown, the criterion of agent $k$ is the infinite horizon, expected total discounted reward, i.e.,  $v^k_{\delta}(\zeta^k(0)) = E\{\sum_t \beta^t r^k(s(t), a(t)) | \zeta^k(0)\}$, where $E\{. |\zeta^k(0)\}$ is the expectation operator conditioned on $\zeta^k(0)$ and where we assume the discount factor $\beta$ is such that $0 \leq \beta < 1$. The discount factor can be interpreted as that there is a probability of $1-\beta$ that the warfare terminates and no further rewards/costs are realized.\\

We assume initial information condition $\zeta(0)$ is given. Let $\Pi^k$ be the set of all finite-memory policies for agent $k, k \in \{L, F_1, F_2\}$, and $\mathscr{Q}$ be the set of all communication matrices $Q=\{P(z(t)|s(t))\}$. We assume each agent wants to maximize the expected total discounted reward. Thus, the goal of follower $i$ is to maximize $v^{F_i}$ by choosing the best response function $\pi^{i,*}: \Pi^L \times \Pi^{F_{3-i}} \times \mathscr{Q} \rightarrow \Pi^{F_i} $ such that $\forall \delta^L \in \Pi^L, Q \in \mathscr{Q}$,
\begin{align}
&v^{F_i}(\delta^{L},\pi^{i,*}(\delta^L, \pi^{3-i,*},Q),\pi^{3-i,*}(\delta^L, \pi^{i,*},Q) )(\zeta^{F_i}(0)) \geq v^{F_i}(\delta^{L},\rho^{F_i}, \pi^{3-i,*}(\delta^L, \rho^{F_i},Q))(\zeta^{F_i}(0)),\notag\\& \forall \rho^{F_i} \in \Pi^{F_i},
i \in \{1,2\}. \label{eq1}
\end{align}

Similarly, the goal of the leader is to maximize $v^L$, taking into consideration the best response policy of each follower; that is,
\begin{align}
&v^L(\pi^L, \pi^{F_1,*}(\pi^L, \pi^{F_2,*},Q),\pi^{F_2,*}(\pi^L, \pi^{F_1,*},Q))(\zeta^{L}(0)) \notag \\&\geq v^L(\delta^L, \pi^{F_1,*}(\delta^L, \pi^{F_2,*},Q),\pi^{F_2,*}(\delta^L, \pi^{F_1,*},Q))(\zeta^{L}(0)), \forall \delta^L \in \Pi^L. \label{eq2}
\end{align}

We remark that searching for the equilibrium policies for the two followers for any given $\delta^L \in \Pi^L$ is a classical decentralized POMDP (Dec-POMDP). Dec-POMDP is NEXP-complete and a survey of solution techniques for Dec-POMDPs can be found in Oliehoek (2012). We focus on the impact of distorted information in communication (as defined by $Q(\epsilon_1, \epsilon_2)$) on  the leader's value function $v^L(\delta^L, \pi^{F_1,*}(\delta^L, \pi^{F_2,*},Q),\pi^{F_2,*}(\delta^L, \pi^{F_1,*},Q))(\zeta^{L}(0))$ for \textit{any} given leader policy $\delta^L \in \Pi^L$ (either optimal or sub-optimal). Developing efficient algorithms to determine optimal leader policies is beyond the scope of this paper and is an interesting research topic.

\section{Value of Distorted Information}
\label{communicationQuality}
In this section, we investigate the impact of distorted information in communication between the two followers on the performance of the leader for any given policy of the leader.
Specifically, we analyze how the leader's value function $v^L$ changes as the probability of having distorted information in communication $\epsilon_1, \epsilon_2$ changes, assuming that a finite-memory leader policy $\delta^L$ and initial conditions $\zeta = (\zeta^L, \zeta^{F_1}, \zeta^{F_2})$ are given. We proceed by first analyzing the best response policies of the two followers.\\

The space of communication quality is $\mathscr{Q}= \{(\epsilon_1,\epsilon_2): 0 \leq \epsilon_1 \leq 1, 0 \leq \epsilon_2\leq 1 \}$. It is clear that the best response policy pair for the followers $(\delta^{F_1}, \delta^{F_2})$ is a function of communication quality $(\epsilon_1, \epsilon_2)$, i.e., $\delta^{F_i} = \pi^{F_i,*}(\delta^L, \pi^{F_{3-i},*},\epsilon_1, \epsilon_2), i \in \{1,2\}$.
\begin{proposition}
	\label{Partition}
	\leavevmode
	Assume $(\delta^{F_1}, \delta^{F_2})$ is the best response policy pair for the followers for communication quality $(\epsilon_1, \epsilon_2)$. Then there exist $c_i \geq 0$ such that if $|\epsilon_i - \epsilon_i'| \leq c_i, i = \{1,2\}$, $(\delta^{F_1}, \delta^{F_2})$ remains the best response policy pair for communication quality $(\epsilon_1', \epsilon_2')$.
\end{proposition}

The values of $c_1,c_2$ is determined in the proof of Proposition \ref{Partition} in Appendix A. Proposition \ref{Partition} indicates that the distorted information will affect the leader's value function depending on the magnitude of distorted information as discussed below.
\begin{enumerate}[(i)]
	\item If the changes of the magnitude of distorted information (in terms of probability $\epsilon_1,\epsilon_2$) are relatively small, the followers will not change their best response policies. Rather, they will select their actions according to the distorted information where the policy of each follower $\delta^{F_i}, i \in \{1,2\}$ remains optimal. A policy $\delta^k:\zeta^{k} \rightarrow A^k$ is a mapping from private information history $\zeta^k$ to the corresponding action space $A^k$. This analysis explores the pure effect of distorted information in communication $\zeta^k$ when the policy of each agent $\delta^k$ is fixed. The leader's performance is affected only by the changes of the followers' actions.
	\item If the changes of the magnitude of distorted information are large enough, the followers may adjust their corresponding best response policies. In this case, at least one of $\delta^{F_i}, i \in \{1,2\}$ is changed. Hence, the distorted information will alter the leader's value function through the adjustment of both policies and actions. 	
\end{enumerate}
As the followers may switch to completely different policies for case (ii), it can be shown by extending the analysis in Chang et al. (2015b) that: (a) the leader's value function may experience discontinuities, and (b) such discontinuities could be favorable or unfavorable to the leader. For example, an attempt by the leader to worsen the quality of communication between adversarial followers may result in the followers switching to a different and potentially less compromised channel of communication, which is unlikely to be to the leader's benefit. As obtaining these results is fairly straightforward, case (ii) is not considered in this paper.\\

Rather, the focus of this paper is on assessing the merit of modulating distorted information in communication for case (i).  We propose and argue that a holistic paradigm, taking into account both the reward structures and the policy employed by each agent is necessary for determining the value of distorted information in communication.  To that end, we employ the structured approach discussed below.

\begin{enumerate}[(1)]
	\item We identify conditions under which the leader's value function $v^L$ can be represented by a power series in $\epsilon_1$ and $\epsilon_2$ (Subsection \ref{sectionPowerSeries}). Thus, the leader's value function is infinitely differentiable. Furthermore, for small $\epsilon_1, \epsilon_2$, the signs of the coefficients of the lower order terms in the power series determine the value of distorted information in communication to the leader.
	\item Assume the conditions of (1) hold:
	\begin{enumerate}[(i)]
		\item 	 We establish the existence of the cases where less accurate information communicated between adversaries (allies) can improve (degrade) the performance of the leader (Proposition \ref{1}) in Subsection \ref{SectionImpact}, by investigating simple zero-memory policies in a totally cooperative/non-cooperative game. More importantly, we demonstrate by example cases that: (a) degraded communication quality between adversaries (allies) does not necessarily lead to an improved (a worsened) leader's performance, and (b) the value of distorted information depends on the specific reward structures among the agents (i.e., zero-sum, general-sum, totally collaborative) and the policy employed by each agent.
		\item We relax the assumptions of zero-memory policies and the perfectly aligned (directly opposite) reward structure in (i) in order to identify a range of reward structures and general finite-memory policies for which the intuitive value of distorted information remains intact. We do so by investigating to what extent the reward structure and policies can deviate from the ones identified in (i).
		\item We also support our thesis that the value of information depends on the policies employed, by demonstrating the simple fact that for any reward structure, the information has no value when the followers totally ignore it in their responses.
		
	\end{enumerate}
	
	\item Finally, in Subsection \ref{othercases}, we discuss the case when the conditions of (1) do not hold, namely, when the power series does not converge or higher order terms in the series have to be considered.
\end{enumerate}

\subsection{The Leader's Value Function as a Power Series in $(\epsilon_1, \epsilon_2)$}
\label{sectionPowerSeries}
It is easy to show by standard POMDP techniques that $\{\zeta^L(t,\tau), y^L(t)\}$ is a sufficient statistic for $\{\zeta^L(t)\}$, where $y^L(t) = P(\zeta^F(t,\tau)|\zeta^L(t))$. Hence, $v^L_{\delta}(\zeta^L(t)) = v^L_{\delta}(\zeta^L(t,\tau), y^L(t))$. We suppress the notation for time $t$ for notational simplicity throughout the rest of the paper. Let $\alpha^{k,l}$ be a real-valued vector over $ \{\zeta\}, \forall k, l \geq 0$. For any vector $u \in R^n$, we say $u \leq 0$, if $u(i) \leq 0, \forall i = 1,...,n$; $u \geq 0$ if $u(i) \geq 0, \forall i = 1,...,n$; and $u$ is unsigned otherwise.

\begin{lemma}
	\label{littleLemmaforConvergence}
	Assume $a^{k,l} \leq c(a^{k-1,l} + a^{k,l-1}) + 2ca^{k-1,l-1}$, where $a,c > 0$, $a^{0,0} \leq 1$, $a^{0,l} \leq c^l$, and $a^{k,0} \leq c^k$.
	\begin{enumerate}[(i)]
		\item For $c \geq 1/4$, $a^{k,l} \leq (4c)^{l+k}$, $\forall k, l \geq 0$;
		\item For $0 < c < 1/4$, $a^{k,l} \leq (4c)^{\max(k,l)}$, $\forall k, l \geq 0$.
	\end{enumerate}
\end{lemma}

\begin{theorem}
	There is a unique real-valued function $g$ on $\{\zeta\}$ satisfying $v^L(\zeta^L, y^L) = \sum_{\zeta^F}y^L(\zeta^F)g^L(\zeta^L,\zeta^F)$. Furthermore, $g^L$ can be represented by a sequence of vectors $\{\alpha^{k,l}\}_{k,l \geq 0}$ as $g^L=\sum_{k=0}^\infty \sum_{l=0}^\infty \epsilon_1^k \epsilon_2^l \alpha^{k,l}$, where the sequence converges for:
	\begin{enumerate}[(i)]
		\item $\forall 0<\epsilon_1,\epsilon_2 <1$, if $0<\beta \leq 1/9$;
		\item $\forall 0 < \epsilon_1, \epsilon_2 < \frac{1-\beta}{8\beta}$,  if $1/9 < \beta < 1$.
	\end{enumerate}
	Hence, $v^L(\zeta^L, y^L)$ is infinitely differentiable in $\epsilon_1, \epsilon_2$.
	\label{thm1}
\end{theorem}

The proof of Theorem \ref{thm1} consists of three major steps. In the first step, we show that the function $g$ is the unique solution to an equation analogous to POMDPs with finite-memory policies. Mathematical induction is then employed to represent $g$ as a power series of $\epsilon_1, \epsilon_2$. The convergence result follows Lemma \ref{littleLemmaforConvergence}. Detailed proof is presented in the Appendix.\\

Theorem \ref{thm1} implies that the leader' value function is infinitely differentiable in $\epsilon_1, \epsilon_2$. Furthermore, for $\epsilon_1, \epsilon_2$ sufficiently small, the value of distorted information in communication can be easily determined by determining the signs of $\alpha^{1,0}$ and $\alpha^{0,1}$. If the joint policies $(\delta^L, \delta^{F_1}, \delta^{F_2})$ are such that $\alpha^{1,0} \leq 0$ (or $\alpha^{0,1} \leq 0$), then the leader's value function can be improved by reducing the probability of having distorted information in communication between two followers, i.e., with smaller $\epsilon_1$ (or $\epsilon_2$). However, if $\alpha^{1,0} \geq 0$ (or $\alpha^{0,1} \geq 0$), then decreasing the probability of having distorted information $\epsilon_1$ (or $\epsilon_2$) may decrease the leader's value function. If $\alpha^{1,0}$ (or $\alpha^{0,1}$) is unsigned, then the value of $y^L$ will determine whether it is worthwhile and, if so, in which direction to change the probability of distorted information communicated between the two followers.

\subsection{The Impact of Reward Structures and Policies}
\label{SectionImpact}
Below, we first demonstrate that there exist reward structures and simple finite memory policies (zero-memory) under which $\alpha^{1,0},\alpha^{0,1} \leq (\geq) 0$. We then discuss the impact of reward structures and policies on the value of distorted information. Let $Q=Q(\epsilon_1, \epsilon_2)$ and $Q'=Q'(\epsilon_1', \epsilon_2')$ where $\epsilon_i \leq \epsilon_i', i = \{1,2\}$.
\subsubsection{The Existence}
\begin{proposition}
	\label{1}
	Assume $r^{F_1} = r^{F_2} = r^F$ and a zero-memory leader policy $\delta^L:S\rightarrow A^L$ is given. Let $\delta^{F_1*}:S^L_1 \times Z^{F_1} \times S^F \rightarrow A^{F_1}$ and $\delta^{F_2*}:S^L_2 \times Z^{F_2} \times S^F \rightarrow A^{F_2}$ be the followers' best response zero-memory policies, i.e., they satisfy
	\begin{align}
	v^{F*}(s) = \max_{a^{F} } \bigg \{R^F(s,\delta^L(s),a^{F}) + \beta \sum_{s'}P(s'|s,\delta^L(s),a^F)v^{F*}(s') \bigg \}, \label{3}
	\end{align}
	For $\epsilon_i, \epsilon_i'$ sufficiently small, $i=1,2$,
	\begin{enumerate}[(1)]
		\item If $r^L = r^F$, then $\alpha^{1,0},\alpha^{0,1} \leq 0$, and $v^L_{Q'}(\zeta^L, y^L) \leq v^L_Q(\zeta^L, y^L)$;
		\item If $r^L= -r^F$, then $\alpha^{1,0}, \alpha^{0,1} \geq 0$, and $v^L_{Q'}(\zeta^L, y^L) \geq v^L_Q(\zeta^L, y^L)$ .
	\end{enumerate}
\end{proposition}

Intuitively, communicating more accurate information between two allies should improve the leader's performance; communicating less accurate information between two adversaries should improve the leader's performance. Proposition \ref{1} verifies this intuition under the conditions that (i) the reward structures are perfectly aligned (or exactly the opposite), and (ii) the followers are using their best response zero-memory policies for $\epsilon_1=\epsilon_2=0$ (the policies remain optimal for small $\epsilon_1,\epsilon_2 > 0$ as guaranteed by Proposition \ref{Partition}). Although this result is not unexpected, Proposition \ref{1} also implies the necessity of taking into account reward structures and policies to avoid misusing distorted information.  Below, we further demonstrate that the intuition may not be always true. Examples 1 and 2 consider simple scenarios depicted in Figure \ref{difference}(b,c) in order to illustrate the importance of the reward structures, while Examples 3 and 4 further analyze the impact of policies. In order to reduce computational complexity, we assume the leader's policy $\delta^L : S \rightarrow A^L$, is fixed. The meaning and the value of each parameter are provided in the Appendix.\\

\begin{example}(\textit{Impact of Reward Structures: Non-cooperative})
	\label{OffensiveExample}
	For a simple case of the scenario in Figure \ref{difference}(b), consider an area held by two groups of militants, each in its own position monitoring army activities in vicinity of the position. The intensity of army activity in vicinity of a position can be high (H) or low (L). Thus, the state space of the army is $S^L = {(H,H), (H,L), (L,H), (L,L)}$. Let $H=1, L=0$ for simplicity.  The militants dynamically move their limited number of heavy weapons between the two positions to meet their needs.  Hence, the state space of the militants is $S^F$ = \{heavy weapons with group 1 (= 0), heavy weapons with group 2 (= 1)\}.  Each militant group collects information on the activity level of the army in vicinity of its position and shares it to the other group. Thus, $Z^i=S^L_i$. The policy of the army is predetermined and given, e.g., units may reconnoiter in vicinity of a position, troops may be moved to the position in the area where the militants have no heavy weapons, etc.  The actions available to a group of militants are to (i) use all weapons including heavy weapons if available (= 0), or (ii) use light weapons and move heavy weapons to the other group (= 1).\\
	
	The reward structure of the army $r^L(s,a)$ consists of two parts: the number of army troops killed by the militants at each decision epoch and a measure of the area reclaimed.  The reward structure of the militants can be difficult to specify as it depends on a wide range of (possibly unknown) issues, e.g., consequence/damage, impact of propaganda and symbolic purposes, etc. (Bier et al. 2007).  We thus consider two games: (i) the reward of the militants is directly opposite to that of the army, i.e., $r^F = -r^L$ (a zero-sum game), and (ii) the signs of $r^L(s,a)$ and $r^F(s,a)$ are the opposite but their absolute values are not the same, i.e., $|r^L| \neq |r^F|$. Namely, a failure to the army is a success to the militants, however, the army places different values on failure compared to the values that the militants place on their own successes.  In each game the followers' employ policies that are optimal (under no distorted information).  The numerical results found in Table \ref{Offensive} show $\alpha^{1,0}, \alpha^{0,1} \geq 0$ if $r^L=-r^F$ as per Proposition \ref{1}. However, if $|r^L|\neq |r^F|$, $\alpha^{1,0}$ and $\alpha^{0,1}$ can be unsigned (neither positive nor negative). Consequently, distorting information shared between adversaries in such situations may not lead to an improvement of the leader's performance. Furthermore, the value of distorted information depends on the value of the leader's belief $y^L$ as well.
\end{example}
\begin{table}
	\begin{center}{\scriptsize
			\begin{tabular}{ | c | c | c | c | c | c | c | c | c | c |}\hline
				\multirow{2}{*}{$\zeta$}&\multicolumn{2}{c|}{Game 1}  & \multicolumn{2}{c|}{Game 2} &  	\multirow{2}{*}{$\zeta$}&\multicolumn{2}{c|}{Game 1}  & \multicolumn{2}{c|}{Game 2} \\
				\cline{2-5}
				\cline{7-10}
				&$\alpha^{0,1}$ & $\alpha^{1,0}$ &  $\alpha^{0,1}$ & $\alpha^{1,0}$ &
				&$\alpha^{0,1}$ & $\alpha^{1,0}$ &  $\alpha^{0,1}$ & $\alpha^{1,0}$\\
				\hline
				[0 0] [0 0] 0	&	10.77	&	5.11	&	2.90	&	0.14	&	[1 0] [0 0] 0	&	13.00	&	6.37	&	-0.50	&	0.15	\\\hline
				[0 1] [0 0] 0	&	11.81	&	3.24	&	2.90	&	0.14	&	[1 1] [0 0] 0	&	13.00	&	6.37	&	-0.50	&	0.15	\\\hline
				[0 0] [1 0] 0	&	12.80	&	3.87	&	2.90	&	0.14	&	[1 0] [1 0] 0	&	10.86	&	4.77	&	-0.50	&	0.15	\\\hline
				[0 1] [1 0] 0	&	14.39	&	5.81	&	2.90	&	0.14	&	[1 1] [1 0] 0	&	10.86	&	4.77	&	-0.50	&	0.15	\\\hline
				[0 0] [0 0] 1	&	9.38	&	3.19	&	-0.54	&	0.02	&	[1 0] [0 0] 1	&	11.19	&	4.02	&	-0.69	&	-0.22	\\\hline
				[0 1] [0 0] 1	&	9.38	&	3.19	&	-3.97	&	0.01	&	[1 1] [0 0] 1	&	11.19	&	4.02	&	-0.69	&	-0.22	\\\hline
				[0 0] [1 0] 1	&	9.38	&	3.19	&	-0.54	&	0.02	&	[1 0] [1 0] 1	&	11.19	&	4.02	&	-0.69	&	-0.22	\\\hline
				[0 1] [1 0] 1	&	9.38	&	3.19	&	-3.97	&	0.01	&	[1 1] [1 0] 1	&	11.19	&	4.02	&	-0.69	&	-0.22	\\\hline
				[0 0] [0 1] 0	&	12.20	&	4.79	&	-1.05	&	0.00	&	[1 0] [0 1] 0	&	13.13	&	4.04	&	-1.93	&	-0.09	\\\hline
				[0 1] [0 1] 0	&	14.70	&	3.70	&	-1.05	&	0.00	&	[1 1] [0 1] 0	&	13.13	&	4.04	&	-1.93	&	-0.09	\\\hline
				[0 0] [1 1] 0	&	12.20	&	4.79	&	1.59	&	0.01	&	[1 0] [1 1] 0	&	13.13	&	4.04	&	-1.75	&	0.15	\\\hline
				[0 1] [1 1] 0	&	14.70	&	3.70	&	1.59	&	0.01	&	[1 1] [1 1] 0	& 	13.13	&	4.04	&	-1.75	&	0.15	\\\hline
				[0 0] [0 1] 1	&	13.66	&	4.19	&	-2.67	&	-0.29	&	[1 0] [0 1] 1	&	15.26	&	3.56	&	-1.87	&	0.08	\\\hline
				[0 1] [0 1] 1	&	13.66	&	4.19	&	1.27	&	0.03	&	[1 1] [0 1] 1	&	15.26	&	3.56	&	-1.87	&	0.08	\\\hline
				[0 0] [1 1] 1	&	13.66	&	4.19	&	0.15	&	-0.22	&	[1 0] [1 1] 1	&	15.26	&	3.56	&	2.90	&	0.02	\\\hline
				[0 1] [1 1] 1	&	13.66	&	4.19	&	-0.44	&	-0.11	&	[1 1] [1 1] 1	&	15.26	&	3.56	&	2.90	&	0.02		
				\\\hline				
		\end{tabular}}
	\end{center}
	\caption{$\alpha$ vectors for the non-cooperative games in Example \ref{OffensiveExample} where $\zeta = ([s^L_1, z^{F_1}], [ z^{F_2},s^L_2], s^F$)}
	\label{Offensive}
\end{table}

\begin{example}(\textit{Impact of Reward Structures: Cooperative})
	\label{DefensiveExample}
	For a simple case of the scenario depicted by Figure \ref{difference}(c), consider a military force consisting of a command element (the leader) and two combat units (allied followers).  Each combat unit is opposed by a group of militants.  The operations of a combat unit against a group of militants can be defined as either close (C) or deep (D).  Formally, deep operations are operations conducted against uncommitted enemy forces, while close operations are operations directed against enemy forces already engaged in battle (Department of the Army 2016).  Here we use the terms close and deep simply to denote one place on the battlefield versus another.  In an urban setting close and deep may be only a few blocks apart making fire support potentially lethal to friendly forces if the information communicated between units about battlefield position is not accurate.  The state space of the command element is $S^L = \{(C, C), (D, C), (C, D), (D, D)\}$. Again, for simplicity, let $C=0, D=1$.  As the militants possess limited heavy weaponry, we assume only one group of militants has heavy weapons at any point in time, hence $S^F$ = \{combat unit 1 facing heavy weapons (= 0), combat unit 2 facing heavy weapons (= 1)\}.  Each combat unit communicates to the other one information about its position and implicitly the state of the command element.  Thus, $Z^i=S^L_i$.  The command element formulates the plan of attack on the militants. The actions available to each combat unit are to attack either (i) the militant group opposing it (= 0), or (ii) the militant group opposing the other combat unit through the provision of fire support (= 1).\\
	
	In warfare it is not uncommon for military forces to be composed of troops from different nations and thus the reward structure of each agent may not be exactly the same. For example, the reward structure of the lead command is a function of terrain taken per period, while the reward structure of each combat unit is a function of terrain taken and the number of casualties incurred.  We thus again consider two games.  In each game the followers' policies are optimal (under no distorted information).  The first game assumes $r^L = r^{F_1}=r^{F_2}$, i.e., a totally cooperative game.  The results in Table \ref{Defensive} show that $\alpha^{1,0}, \alpha^{0,1} \le 0$, as dictated by Proposition \ref{1}. The second game is where the signs of $r^L(s,a)$ and $r^F(s,a)$ are the same but $r^L \neq r^F$.  The numerical results show that $\alpha^{1,0}$ and $\alpha^{0,1}$ can be unsigned. Consequently, more accurate information sharing between the two allies does not guarantee an improvement in the leader's performance.
\end{example}

\begin{table}
	\begin{center}{\scriptsize
			\begin{tabular}{ | c | c | c | c | c | c | c | c | c | c |}\hline
				\multirow{2}{*}{$\zeta$}&\multicolumn{2}{c|}{Game 1}  & \multicolumn{2}{c|}{Game 2} &  	\multirow{2}{*}{$\zeta$}&\multicolumn{2}{c|}{Game 1}  & \multicolumn{2}{c|}{Game 2} \\
				\cline{2-5}
				\cline{7-10}
				&$\alpha^{0,1}$ & $\alpha^{1,0}$ &  $\alpha^{0,1}$ & $\alpha^{1,0}$ &
				&$\alpha^{0,1}$ & $\alpha^{1,0}$ &  $\alpha^{0,1}$ & $\alpha^{1,0}$\\
				\hline
				[0 0] [0 0] 0	&	-20.49	&	-7.57	&	-2.74	&	-2.81	&	[1 0] [0 0] 0	&	-21.66	&	-8.18	&	2.10	&	0.93	\\\hline
				[0 1] [0 0] 0	&	-20.49	&	-7.57	&	-2.74	&	-2.81	&	[1 1] [0 0] 0	&	-20.79	&	-6.82	&	1.45	&	0.08	\\\hline
				[0 0] [1 0] 0	&	-22.55	&	-10.37	&	2.35	&	2.11	&	[1 0] [1 0] 0	&	-21.03	&	-7.97	&	1.35	&	1.83	\\\hline
				[0 1] [1 0] 0	&	-22.55	&	-10.37	&	2.35	&	2.11	&	[1 1] [1 0] 0	&	-18.99	&	-8.73	&	-1.49	&	-1.14	\\\hline
				[0 0] [0 0] 1	&	-23.96	&	-8.85	&	-2.00	&	-1.91	&	[1 0] [0 0] 1	&	-20.86	&	-8.37	&	1.84	&	1.21	\\\hline
				[0 1] [0 0] 1	&	-23.96	&	-8.85	&	-2.00	&	-1.91	&	[1 1] [0 0] 1	&	-20.07	&	-6.49	&	1.84	&	1.21	\\\hline
				[0 0] [1 0] 1	&	-22.42	&	-9.63	&	-2.00	&	-1.91	&	[1 0] [1 0] 1	&	-17.68	&	-6.8	&	1.84	&	1.21	\\\hline
				[0 1] [1 0] 1	&	-22.42	&	-9.63	&	-2.00	&	-1.91	&	[1 1] [1 0] 1	&	-18.79	&	-6.34	&	1.84	&	1.21	\\\hline
				[0 0] [0 1] 0	&	-15.97	&	-7.67	&	1.31	&	0.50	&	[1 0] [0 1] 0	&	-21.25	&	-6.48	&	3.78	&	3.11	\\\hline
				[0 1] [0 1] 0	&	-15.97	&	-7.67	&	1.31	&	0.50	&	[1 1] [0 1] 0	&	-21.93	&	-8.11	&	3.93	&	3.50	\\\hline
				[0 0] [1 1] 0	&	-23.98	&	-9.24	&	0.63	&	0.55	&	[1 0] [1 1] 0	&	-17.95	&	-7.20	&	0.79	&	-0.52	\\\hline
				[0 1] [1 1] 0	&	-23.98	&	-9.24	&	0.63	&	0.55	&	[1 1] [1 1] 0	& 	-18.53	&	-10.25	&	-2.39	&	-3.43	\\\hline
				[0 0] [0 1] 1	&	-22.37	&	-9.04	&	2.99	&	2.00	&	[1 0] [0 1] 1	&	-21.79	&	-7.78	&	0.15	&	-0.86	\\\hline
				[0 1] [0 1] 1	&	-22.37	&	-9.04	&	2.99	&	2.00	&	[1 1] [0 1] 1	&	-19.9	&	-4.99	&	0.15	&	-0.86	\\\hline
				[0 0] [1 1] 1	&	-22.37	&	-9.04	&	0.68	&	0.08	&	[1 0] [1 1] 1	&	-21.79	&	-7.78	&	1.24	&	1.28	\\\hline
				[0 1] [1 1] 1	&	-22.37	&	-9.04	&	0.68	&	0.08	&	[1 1] [1 1] 1	&	-19.9	&	-4.99	&	1.24	&	1.28		
				\\\hline				
		\end{tabular}}
	\end{center}
	\caption{$\alpha$ vectors for the cooperative games in Example \ref{DefensiveExample} where $\zeta = ([s^L_1, z^{F_1}], [ z^{F_2},s^L_2], s^F$)}
	\label{Defensive}
\end{table}
\begin{example}
	\label{MyopicExample}
	(\textit{Impact of Policies: Myopic Policy})
	Example \ref{DefensiveExample} continued. We now compare the values of $\alpha^{1,0}$ and $\alpha^{0,1}$ under two different followers' response policies in Table \ref{myopictable}, assuming $r^L=r^{F_1}=r^{F_2}$. The first response policies satisfy the optimality equation \eqref{3}. Hence, all elements of $\alpha^{0,1}$ and $\alpha^{1,0}$ are non-positive as per Proposition \ref{1}, indicating more accurate information in communication between two allies can indeed improve the leader's performance. The second response policies are myopic in that the two followers will select the action $a^{F*}$ such that $r^F(s,\delta^{L}(s),a^{F*})=\max_{a^F \in A^F} r^F(s,\delta^{L}(s),a^F)$ for each $s, \forall t$. The signs of $\alpha^{0,1}$ and $\alpha^{1,0}$ are, however, unsigned. The value of distorted information will also depend on the belief vector $y^L$.
\end{example}

\begin{table}
	\begin{center}
		{\scriptsize
			\begin{tabular}{ | c | c | c | c | c | c | c | c | c | c |} \hline			
				\multirow{2}{*}{$\zeta$} & \multicolumn{2}{c|}{Optimal in \eqref{3}} &\multicolumn{2}{c|}{Myopic} &  	\multirow{2}{*}{$\zeta$}  & \multicolumn{2}{c|}{Optimal in \eqref{3}} &\multicolumn{2}{c|}{Myopic}\\
				\cline{2-5}
				\cline{7-10}
				&$\alpha^{0,1}$ & $\alpha^{1,0}$ &  $\alpha^{0,1}$ & $\alpha^{1,0}$ &
				&$\alpha^{0,1}$ & $\alpha^{1,0}$ &  $\alpha^{0,1}$ & $\alpha^{1,0}$\\
				\hline
				[0 0] [0 0] 0	&	-19.16	&	-4.85	&	0.62	&	0.62	&	[1 0] [0 0] 0	&	-13.11	&	-5.59	&	-0.82	&	-0.74	\\\hline
				[0 1] [0 0] 0	&	-19.16	&	-4.85	&	0.62	&	0.62	&	[1 1] [0 0] 0	&	-13.11	&	-5.59	&	-0.82	&	-0.74	\\\hline
				[0 0] [1 0] 0	&	-14.05	&	-4.65	&	0.62	&	0.62	&	[1 0] [1 0] 0	&	-13.87	&	-5.19	&	-0.82	&	-0.74	\\\hline
				[0 1] [1 0] 0	&	-14.05	&	-4.65	&	0.62	&	0.62	&	[1 1] [1 0] 0	&	-13.87	&	-5.19	&	-0.82	&	-0.74	\\\hline
				[0 0] [0 0] 1	&	-15.44	&	-5.23	&	0.03	&	0.03	&	[1 0] [0 0] 1	&	-10.72	&	-5.55	&	0.24	&	0.24	\\\hline
				[0 1] [0 0] 1	&	-8.64	&	-7.55	&	-0.76	&	-0.76	&	[1 1] [0 0] 1	&	-10.72	&	-5.55	&	0.24	&	0.24	\\\hline
				[0 0] [1 0] 1	&	-15.44	&	-5.23	&	-0.76	&	-0.76	&	[1 0] [1 0] 1	&	-10.72	&	-5.55	&	0.81	&	0.81	\\\hline
				[0 1] [1 0] 1	&	-8.64	&	-7.55	&	1.00	&	1.00	&	[1 1] [1 0] 1	&	-10.72	&	-5.55	&	0.81	&	0.81	\\\hline
				[0 0] [0 1] 0	&	-13.87	&	-5.19	&	-0.74	&	-0.82	&	[1 0] [0 1] 0	&	-23.5	&	-3.34	&	-0.52	&	-0.52	\\\hline
				[0 1] [0 1] 0	&	-13.87	&	-5.19	&	-0.74	&	-0.82	&	[1 1] [0 1] 0	&	-23.5	&	-3.34	&	-0.52	&	-0.52	\\\hline
				[0 0] [1 1] 0	&	-13.11	&	-5.69	&	-0.74	&	-0.82	&	[1 0] [1 1] 0	&	-17.46	&	-6.34	&	-0.52	&	-0.52	\\\hline
				[0 1] [1 1] 0	&	-13.11	&	-5.69	&	-0.74	&	-0.82	&	[1 1] [1 1] 0	&	-17.46	&	-6.34	&	-0.52	&	-0.52	\\\hline
				[0 0] [0 1] 1	&	-16.5	&	-4.27	&	0.24	&	0.24	&	[1 0] [0 1] 1	&	-12.9	&	-5.53	&	-0.57	&	-0.57	\\\hline
				[0 1] [0 1] 1	&	-10.72	&	-5.55	&	0.81	&	0.81	&	[1 1] [0 1] 1	&	-12.9	&	-5.53	&	-0.57	&	-0.57	\\\hline
				[0 0] [1 1] 1	&	-16.5	&	-4.27	&	0.24	&	0.24	&	[1 0] [1 1] 1	&	-12.9	&	-5.53	&	-0.57	&	-0.57	\\\hline
				[0 1] [1 1] 1	&	-10.72	&	-5.55	&	0.81	&	0.81	&	[1 1] [1 1] 1	&	-12.9	&	-5.53	&	-0.57	&	-0.57	
				\\\hline				
		\end{tabular}}	
	\end{center}
	\caption{$\alpha$ vector for a myopic policy in Example \ref{MyopicExample} where $\zeta = ([s^L_1, z^{F_1}], [ z^{F_2},s^L_2], s^F$)}
	\label{myopictable}
\end{table}
\begin{example}
	\label{kstepEx}
	(\textit{Impact of Policies: $k$-Step Ahead Policy})
	Example \ref{DefensiveExample} continued. Assuming $r^L=r^{F_1}=r^{F_2}$, we further compare the values of $\alpha^{1,0}$ and $\alpha^{0,1}$ under a $k$-step ahead policy in Table \ref{ksteptable}. The optimal policy (under no distorted information) is also a 4-step ahead policy, where all elements of $\alpha^{0,1}$ and $\alpha^{1,0}$ are negative as expected. It is noteworthy that even under a 3-step ahead policy (just 1-step away from optimal), the vector $\alpha^{1,0}$ is still greater than zero.  Thus, more accurate information sharing between allies may not necessarily improve the leader's value function even under some good sub-optimal policies.
\end{example}

\begin{table}
	\begin{center}
		{\scriptsize
			\begin{tabular}{ | c | c | c | c | c | c | c | c | c | c |} \hline			
				\multirow{2}{*}{$\zeta$} & \multicolumn{2}{c|}{Optimal in \eqref{3}} &\multicolumn{2}{c|}{3-step ahead} &  	\multirow{2}{*}{$\zeta$}  & \multicolumn{2}{c|}{Optimal in \eqref{3}} &\multicolumn{2}{c|}{3-step ahead}\\
				\cline{2-5}
				\cline{7-10}
				&$\alpha^{0,1}$ & $\alpha^{1,0}$ &  $\alpha^{0,1}$ & $\alpha^{1,0}$ &
				&$\alpha^{0,1}$ & $\alpha^{1,0}$ &  $\alpha^{0,1}$ & $\alpha^{1,0}$\\
				\hline
				[0 0] [0 0] 0	&	-23736.86	&	-406.95	&	-17492.31	&	211.92	&	[1 0] [0 0] 0	&	-23612.46	&	-407.57	&	-17286.38	&	220.08	\\\hline
				[0 1] [0 0] 0	&	-23736.86	&	-406.95	&	-17492.31	&	211.92	&	[1 1] [0 0] 0	&	-23612.46	&	-407.57	&	-17286.38	&	220.08	\\\hline
				[0 0] [1 0] 0	&	-23736.86	&	-406.95	&	-17492.31	&	211.92	&	[1 0] [1 0] 0	&	-23612.46	&	-407.57	&	-17286.38	&	220.08	\\\hline
				[0 1] [1 0] 0	&	-23736.86	&	-406.95	&	-17492.31	&	211.92	&	[1 1] [1 0] 0	&	-23612.46	&	-407.57	&	-17286.38	&	220.08	\\\hline
				[0 0] [0 0] 1	&	-23727.83	&	-404.17	&	-17471.76	&	217.15	&	[1 0] [0 0] 1	&	-23722.75	&	-406.78	&	-17463.76	&	214.73	\\\hline
				[0 1] [0 0] 1	&	-23622.41	&	-407.81	&	-17293.29	&	221.15	&	[1 1] [0 0] 1	&	-23722.75	&	-406.78	&	-17500.42	&	213.6	\\\hline
				[0 0] [1 0] 1	&	-23727.83	&	-404.17	&	-17471.76	&	217.15	&	[1 0] [1 0] 1	&	-23722.75	&	-406.78	&	-17463.76	&	214.73	\\\hline
				[0 1] [1 0] 1	&	-23622.41	&	-407.81	&	-17293.29	&	221.15	&	[1 1] [1 0] 1	&	-23722.75	&	-406.78	&	-17500.42	&	213.6	\\\hline
				[0 0] [0 1] 0	&	-23612.46	&	-407.57	&	-17286.38	&	220.08	&	[1 0] [0 1] 0	&	-23665.64	&	-409.74	&	-17413.94	&	208.88	\\\hline
				[0 1] [0 1] 0	&	-23612.46	&	-407.57	&	-17286.38	&	220.08	&	[1 1] [0 1] 0	&	-23665.64	&	-409.74	&	-17413.94	&	208.88	\\\hline
				[0 0] [1 1] 0	&	-23722.95	&	-406.5	&	-17503.54	&	210.87	&	[1 0] [1 1] 0	&	-23883.22	&	-403.81	&	-17705.59	&	209.86	\\\hline
				[0 1] [1 1] 0	&	-23722.95	&	-406.5	&	-17503.54	&	210.87	&	[1 1] [1 1] 0	&	-23883.22	&	-403.81	&	-17705.59	&	209.86	\\\hline
				[0 0] [0 1] 1	&	-23741.84	&	-405.08	&	-17500.42	&	213.6	&	[1 0] [0 1] 1	&	-23795.26	&	-405.59	&	-17604.43	&	207.35	\\\hline
				[0 1] [0 1] 1	&	-23722.75	&	-406.78	&	-17463.76	&	214.73	&	[1 1] [0 1] 1	&	-23795.26	&	-405.59	&	-17232.05	&	224.17	\\\hline
				[0 0] [1 1] 1	&	-23741.84	&	-405.08	&	-17500.42	&	213.6	&	[1 0] [1 1] 1	&	-23795.26	&	-405.59	&	-17604.43	&	207.35	\\\hline
				[0 1] [1 1] 1	&	-23722.75	&	-406.78	&	-17463.76	&	214.73	&	[1 1] [1 1] 1	&	-23795.26	&	-405.59	&	-17232.05	&	224.17	\\\hline			
		\end{tabular}}	
	\end{center}
	\caption{$\alpha$ vector for a $k$-step ahead policy in Example \ref{kstepEx} where $\zeta = ([s^L_1, z^{F_1}], [ z^{F_2},s^L_2], s^F$)}
	\label{ksteptable}
\end{table}

So far, we have employed the assumption of zero-memory policies and established conditions under which the value of distorted information is intuitively sound. If these conditions are violated, however, the intuition may not hold as demonstrated by relevant examples. We further remark these examples are not rare. Out of 500,000 randomly generated examples similar to Examples \ref{OffensiveExample}-\ref{kstepEx} with a variety of sizes for the state spaces, action spaces, and observation spaces, 10.1\% of the examples violated the intuition.

\subsubsection{General Cases}
We now relax the assumption of zero-memory policies and investigate a range of reward structures and general \textit{finite memory} policies within which our intuitive understanding for distorted information still holds. The objective of this analysis to address the question: to what degree can the reward structures and policies deviate from the ones identified in Proposition \ref{1} such that the prevalent intuition remains valid.\\

Assume there exists a reward structure $r^*=(r^{L,*},r^{F_1,*},r^{F_2,*})$ and a policy tuple $\delta^*=(\delta^{L}, \delta^{F_1,*}, \delta^{F_2,*}) \in \Pi=\Pi^L \times \Pi^{F_1} \times \Pi^{F_2}$, where $\alpha^{1,0}_{\delta^*}, \alpha^{0,1}_{\delta^*} < 0$, which we call the ``nominal" reward structure and policy tuple (their existence is guaranteed as per Proposition \ref{1}). We now investigate how far a reward structure $r=(r^{L},r^{F_1},r^{F_2})$ and a policy tuple $\delta=(\delta^L, \delta^{F_1}, \delta^{F_2}) \in \Pi$ can deviate from the nominal ones such that $\alpha^{1,0}_{\delta},\alpha^{0,1}_{\delta} \leq 0$. The conditions for which $\alpha^{1,0}_{\delta},\alpha^{0,1}_{\delta} \geq 0$ can be found in a similar manner and thus are omitted.\\

Let $U$ be the set of all bounded, real-valued functions on $\{\zeta \}$. Let $||u||=\max_\zeta|u(\zeta)|$ for a vector $u$.
$\forall \delta,\delta' \in \Pi$, define the distance function on $\Pi$ by:
\begin{align}
&d(\delta, \delta') \notag\\
&= \sup_{||u|| \leq 1, u \in U}||\sum_{s'}P(s'|s,\delta(\zeta))u(\bar{\zeta}(s',s^{L'}_2,s^{L'}_1,\zeta,\delta(\zeta)))
-\sum_{s'}P(s'|s,\delta'(\zeta))u(\bar{\zeta}(s',s^{L'}_2,s^{L'}_1,\zeta,\delta'(\zeta)))||.
\end{align}

\begin{lemma}
	\label{distanceProposition}
	\begin{enumerate}[(a)]
		\item[]
		\item $||P(\delta) - P(\delta')|| \leq d(\delta, \delta') \leq 2$,
		\item If $\delta$ and $\delta'$ are all zero memory policies, $d(\delta,\delta')= ||P(\delta) - P(\delta')||.$
	\end{enumerate}
\end{lemma}

Assume $\alpha^{1,0}_{\delta^*} < 0$ and $\alpha^{0,1}_{\delta^*} < 0$, and let $h_1=|\max_{\zeta} \alpha^{1,0}_{\delta^*}(\zeta)|$ and $h_2 = |\max_{\zeta} \alpha^{0,1}_{\delta^*}(\zeta)|$.
\begin{theorem}
	\label{deviation}
	Let $\eta^1=||\alpha^{1,0}_{\delta^*}||+||\sum_{z^{F_1} \neq s^{L}_2}\sigma_{s^{L}_2 ,z^{F_1}}\alpha^{0,0}_{\delta^*}(\bar{\zeta}(s,z^{F_1},s^{L}_1,\zeta,\delta^*(\zeta)))-\alpha^{0,0}_{\delta^*}(\bar{\zeta}(s,s^{L}_2,s^{L}_1,\zeta,\delta^*(\zeta)))||$, $\eta^2=||\alpha^{0,1}_{\delta^*}||+||\sum_{z^{F_2}\neq s^{L}_1}\sigma_{s^{L}_1,z^{F_2}}\alpha^{0,0}_{\delta^*}(\bar{\zeta}(s,s^{L}_2, z^{F_2}, \zeta, \delta^*(\zeta))) - \alpha^{0,0}_{\delta^*}(\bar{\zeta}(s,s^{L}_2,s^{L}_1, \zeta, \delta^*(\zeta)))||$.
	Then, for any finite memory policy tuple $\delta \in \Pi$, $\alpha^{1,0}_\delta, \alpha^{0,1}_\delta \leq 0$ if the following two conditions are satisfied:
	\begin{align*}
	&d(\delta,\delta^*) \eta^1 + \left (||\sum_{z^{F_1}\neq s^{L}_2}\sigma_{s^{L}_2,z^{F_1}} \left[\alpha^{0,0}_{\delta^*}(\bar{\zeta}(s, z^{F_1},s^{L}_1, \zeta, \delta(\zeta))) - \alpha^{0,0}_\delta(\bar{\zeta}(s, z^{F_1}, s^{L}_1,\zeta, \delta(\zeta)))\right]||\right ) \\
	&+ \left (||\alpha^{0,0}_{\delta^*}(\bar{\zeta}(s,s^{L}_2,s^{L}_1, \zeta, \delta(\zeta)))- \alpha^{0,0}_{\delta}(\bar{\zeta}(s,s^{L}_2,s^{L}_1, \zeta, \delta(\zeta)))||\right)	\leq \frac{1-\beta }{\beta}h_1,\\
	&d(\delta,\delta^*)\eta^2 + \left (||\sum_{z^{F_2}\neq s^{L}_1}\sigma_{s^{L}_1,z^{F_2}}\left[\alpha^{0,0}_{\delta^*}(\bar{\zeta}(s,s^{L}_2, z^{F_2}, \zeta, \delta(\zeta))) -\alpha^{0,0}_\delta(\bar{\zeta}(s, s^{L}_1,z^{F_2}, \zeta, \delta(\zeta)))\right]||\right ) \notag\\
	&+ \left (||\alpha^{0,0}_{\delta^*}(\bar{\zeta}(s,s^{L}_2,s^{L}_1, \zeta, \delta(\zeta)))- \alpha^{0,0}_{\delta}(\bar{\zeta}(s,s^{L}_2,s^{L}_1, \zeta, \delta(\zeta)))||\right)	\leq \frac{1-\beta }{\beta}h_2.
	\end{align*}
\end{theorem}
\begin{corollary}
	\label{deviationcorollary}
	Assume both $\delta^*$ and $\delta$ are zero-memory policy tuples.\\
	Let $\eta^1=||\alpha^{1,0}_{\delta^*}||+||\sum_{z^{F_1}\neq s^{L}_2}\sigma_{s^{L}_2,z^{F_1}}\alpha^{0,0}_{\delta^*}(\{s,z^{F_1},s^L_1\}) - \alpha^{0,0}_{\delta^*}(\{s,s^L_2,s^L_1\})||$, and  $\eta^2=||\alpha^{0,1}_{\delta^*}||+||\sum_{z^{F_2}\neq s^{L}_1}\sigma_{s^{L}_1,z^{F_2}}\alpha^{0,0}_{\delta^*}(\{s,s^L_2,z^{F_2}\}) - \alpha^{0,0}_{\delta^*}(\{s,s^L_2,s^L_1\}))||$. Then, $\alpha^{1,0}_\delta$, $\alpha^{0,1}_\delta \leq 0$ if the following two conditions are satisfied:
	\begin{align}
	\eta^i ||P(\delta)-P(\delta^*)||+2||\alpha^{0,0}_{\delta^*} - \alpha^{0,0}_\delta|| \leq \frac{1-\beta}{\beta}h_i, i =\{1,2\}.
	\end{align}
\end{corollary}

Theorem \ref{deviation} provides sufficient conditions for any finite memory policies, while Corollary \ref{deviationcorollary} focuses specifically on zero-memory policies. Thus, the value of distorted information for an arbitrary reward structure and a policy tuple can be determined by examining $d(\delta,\delta^*)$ and the distance of $\alpha^{0,0}_\delta$ to the nominal value $\alpha^{0,0}_{\delta^*}$. We remark that the conditions established in both Theorem \ref{deviation} and Corollary \ref{deviationcorollary} are only sufficient conditions and therefore, sometimes can be conservative. Investigating the existence of necessary and sufficient conditions appears to be an interesting research topic. Below, we employ Examples \ref{deviationExampleReward} and \ref{deviationExample} for assessing holistically the impact of reward structure and policy on the value of distorted information.  Specifically,  we illustrate how the value of distorted information changes as we deviate from the conditions established in Proposition \ref{1}.\\

\begin{figure}
	\centering
	\begin{minipage}[b]{0.75\textwidth}
		\includegraphics[width=1\linewidth]{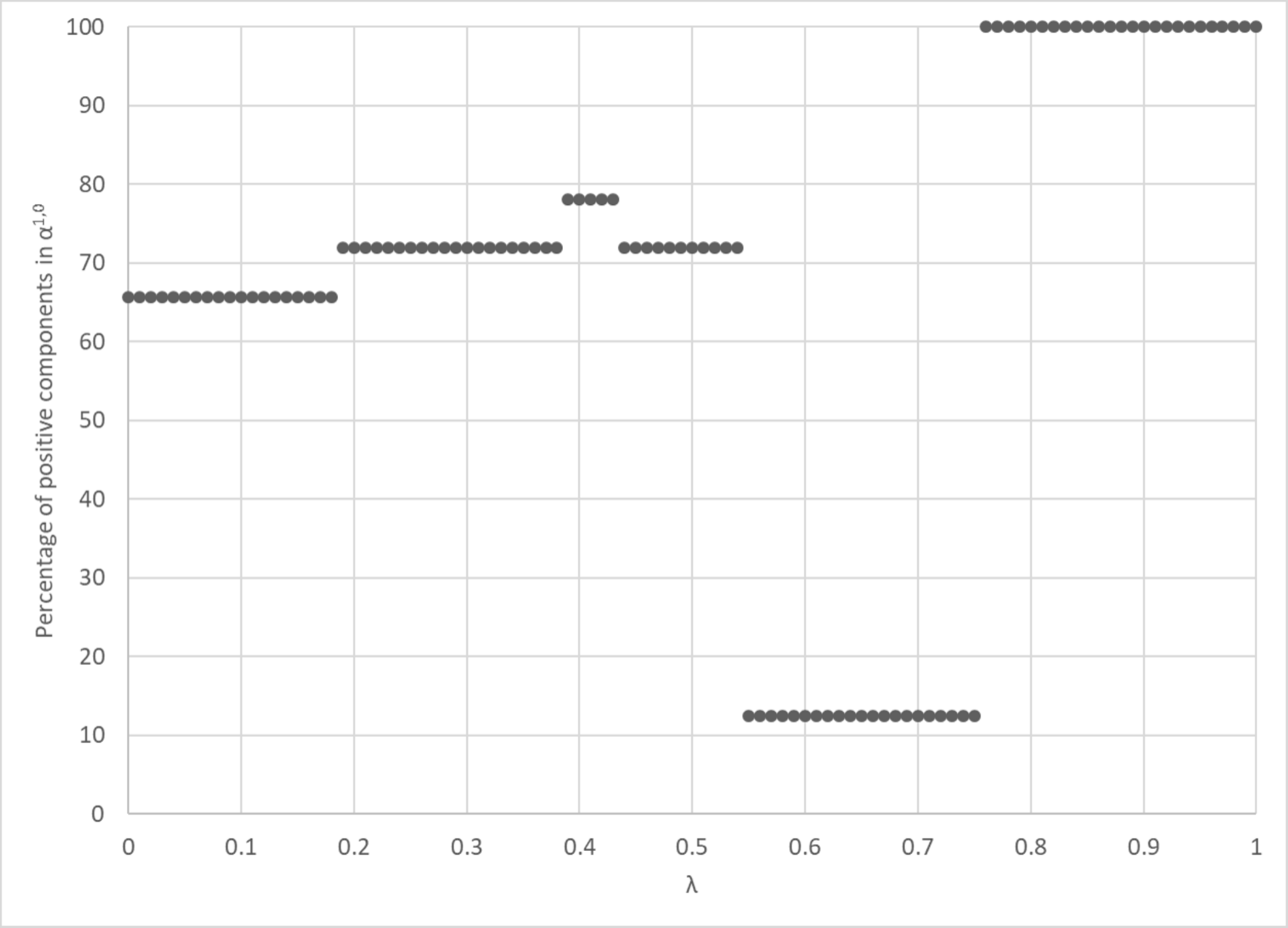}\\
		{(a) non-cooperative situations}
	\end{minipage}
	
	\begin{minipage}[b]{0.75\textwidth}
		\includegraphics[width=1\linewidth]{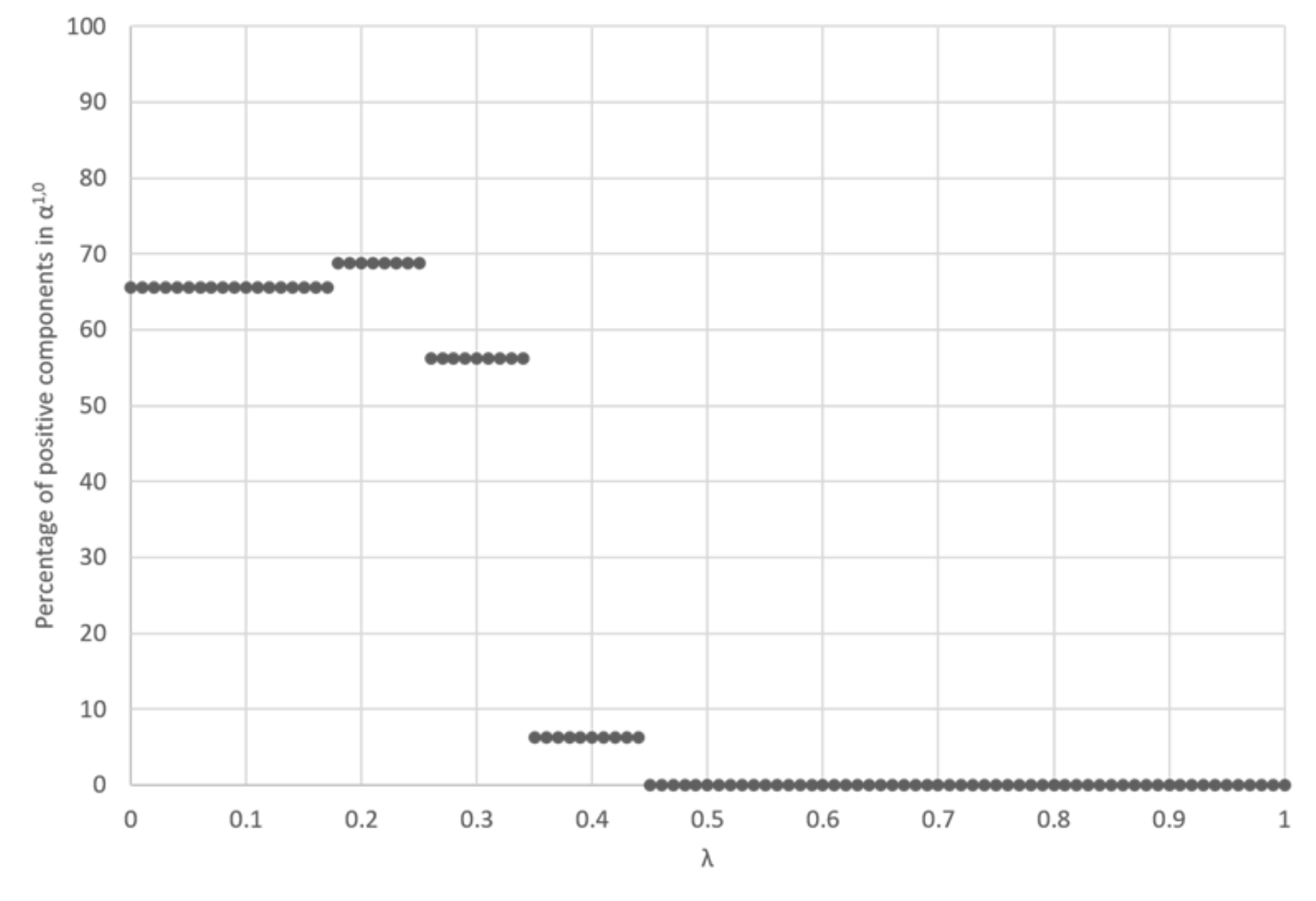}
		\\{(b) cooperative situations}
	\end{minipage}
	\caption{Value of distorted information as the reward structure changes in Example \ref{deviationExampleReward}}
	
	\label{deviationExampleRewardFigure}
\end{figure}

\begin{figure}
	\begin{center}
		\includegraphics[width=0.85\textwidth]{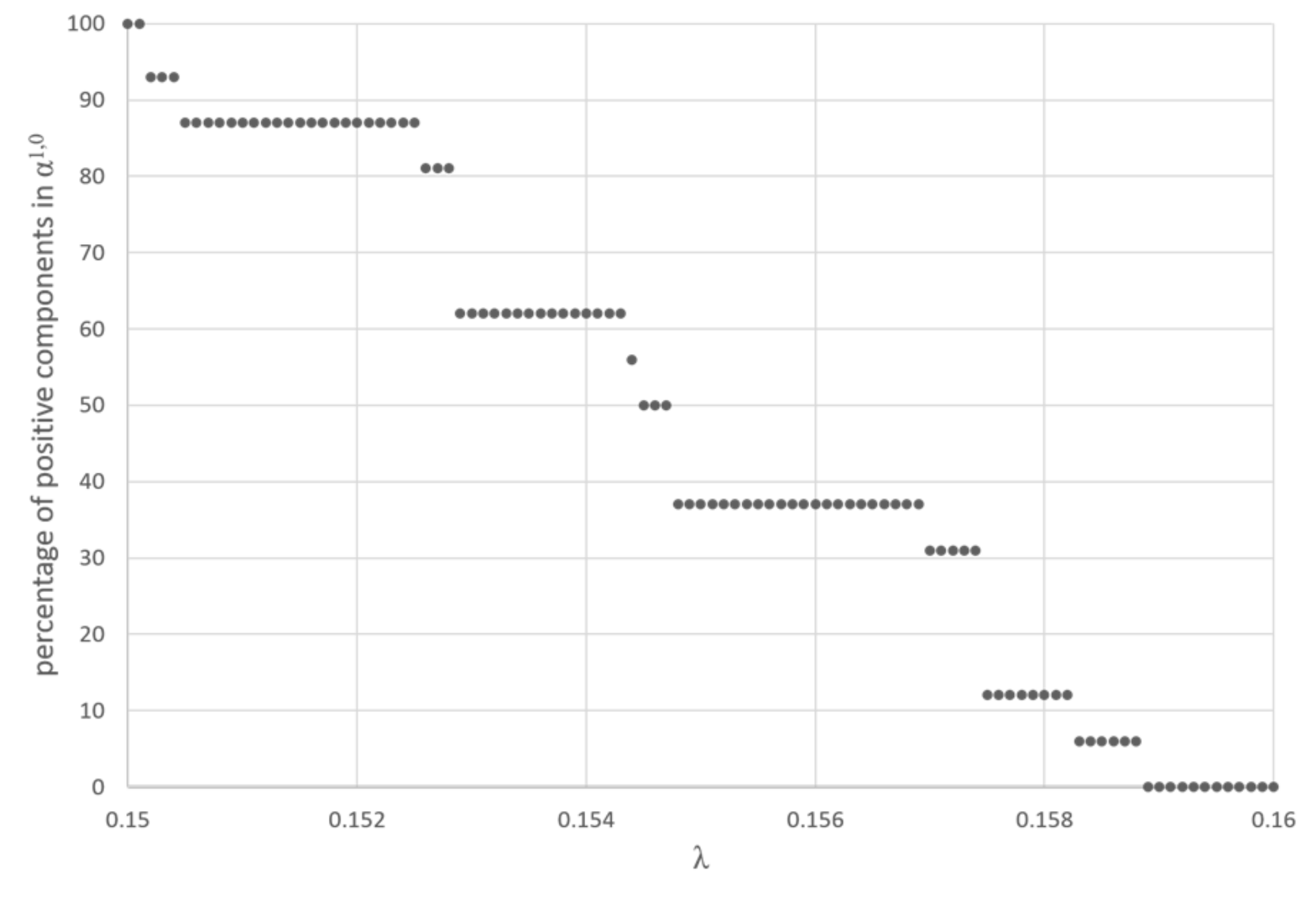}
		\caption{Value of distorted information as policy changes in Example \ref{deviationExample}}
		\label{randomPolicy}
	\end{center}
\end{figure}

\begin{example}
	\label{deviationExampleReward}
	Continuation of Examples \ref{OffensiveExample} and \ref{Defensive}. Let $r^{F,*}=-r^L$ ($r^{F,*}=r^L$) and let $r^{F,nc}$ ($r^{F,c}$) be the reward structure of the two followers in the second game in Example \ref{OffensiveExample} (Example \ref{Defensive}). Construct a series of reward structures for the two followers by $r^F = \lambda r^{F,*} + (1-\lambda) r^{F,nc}$ ($r^F = \lambda r^{F,*} + (1-\lambda) r^{F,c}$), $\lambda \in [0,1]$. For any $\lambda$, optimal response policies are employed by the two followers. Figure \ref{deviationExampleRewardFigure} shows how the value of distorted information evolves as the reward structure changes. Specifically, it shows how the percentage of positive components in $\alpha^{1,0}$ changes as $\lambda$ increases. Note that any game with $\lambda \in [0,1]$ is a non-cooperative example of Figure \ref{difference}(b) (a cooperative example of Figure \ref{difference}(c)); however, $\alpha^{1,0} \geq 0$ ($\alpha^{1,0} \leq 0$) only if the reward structure is close enough to the zero-sum (totally cooperative) case.
\end{example}
\begin{example}
	\label{deviationExample}
	As shown in Example \ref{kstepEx}, for a given leader's policy $\delta^L$, there may not exist a $k$-step ahead policy $\delta^F$ around the optimal policy $\delta^{F,*}$ such that $\alpha^{1,0}_{(\delta^L,\delta^F)} < 0$ and $\alpha^{0,1}_{(\delta^L,\delta^F)} < 0$. We thus construct a series of randomized policies by $\rho^F = (1-\lambda)\delta^{F,3}+\lambda \delta^{F,*}, \lambda \in [0,1]$, where $\delta^{F,3}$ is the 3-step ahead policy pair and  $\delta^{F,*}$ is the optimal policy pair for the followers of Example \ref{kstepEx}. Results in Table \ref{ksteptable} show $\alpha^{1,0}_{(\delta^L, \delta^{F,3})} > 0$ but $\alpha^{1,0}_{(\delta^L, \delta^{F,*})} < 0$ (both $\alpha^{1,0}_{(\delta^L, \delta^{F,3})}$ and $\alpha^{1,0}_{(\delta^L, \delta^{F,*})} \geq 0$). Figure \ref{randomPolicy} shows how the value of distorted information evolves as the policy changes from sub-optimal to optimal for the case under consideration. Specifically, it depicts the percentage of positive components in $\alpha^{1,0}_{(\delta^L, \rho^F)}$ decreasing as $\lambda$ increases (i.e., $\delta^F$ moves from $\delta^{F,3}$ to $\delta^{F,*}$). Clearly, $\alpha^{1,0}_{(\delta^L, \rho^F)} < 0$ only if $\delta^F$ is close enough to $\delta^{F,*}$.
\end{example}
\subsubsection{Zero Value of Distorted Information}
We highlight another implication of Theorem \ref{thm1}: if the followers' policies are independent of real-time information, then the communication of distorted information between followers has no impact on the leader's performance. For example, if adversaries are committed to taking certain actions regardless of the real-time information received, then there is no value to the leader in distorting the information passing between the followers, regardless of the leader's policy $\delta^L$.  This reinforces the idea that policies and the value of distorted information need to be considered and analyzed together.  Proposition \ref{2} captures this conjecture.

\begin{proposition}
	\label{2}
	Assume $\delta^{F_1}(\zeta^{F_1}) = a^{F_1}$ and $\delta^{F_2}(\zeta^{F_2}) = a^{F_2}$ for all $\zeta^{F_1}$ and $\zeta^{F_2}$. Then $\alpha^{k,0} = 0$, $\alpha^{0,l} = 0$, $\alpha^{k,l} = 0$ for $k,l \geq 1$ and $v^L(\zeta^L, y^L) = \sum_{\zeta^F} y^L(\zeta^F) \alpha^{0,0}(\zeta^L,\zeta^F)$. Hence, the value of distorted information is zero for any finite-memory leader policy $\delta^L \in \Pi^L$ and reward structure $r=(r^L, r^{F_1}, r^{F_2})$.
\end{proposition}
\subsection{Other Cases for the Leader's Value Function}
\label{othercases}
For large $\epsilon_1, \epsilon_2$, (i) the infinite series may not converge, especially for a large discount factor $\beta$, or (ii) the linear approximation may not be a good approximation of $g^L$. The value of distorted information may also depend on the signs of higher order $\alpha^{k,l}$s that may not be easy to analyze. Thus, we present another result guaranteeing that more (less) accurate information sharing between two allies (adversaries) will improve the leader's value function.
\begin{theorem}
	\label{convex} Assume:
	\begin{enumerate}[(i)]
		\item there is a stochastic matrix $R^i$ such that $Q^iR^i=Q'^i$;
		\item $v^{F_i}_{\delta,Q}(\zeta^{F_i}(t,\tau), y^{F_i}(t))$ is concave in $y^{F_i}(t)$ for $\zeta^{F_i}(t,\tau)$;
		\item $\zeta^k(0) = \{\zeta^k(0,\tau),y^k(0)\}, k \in \{L,F_1,F_2\}$ and $P(\zeta^{L,F_j}(0)|\zeta^{F_i}(0))=1, i \neq j, i, j \in \{1,2\}$.
	\end{enumerate}
	Then,
	$v^{F_i}_{\delta,Q'}(\zeta^{F_i}(t,\tau), y^{F_i}(t)) \leq v^{F_i}_{\delta,Q}(\zeta^{F_i}(t,\tau), y^{F_i}(t)), i \in \{1,2\}$, hence $v^L_{\delta,Q'}(\zeta^L(0)) \leq v^L_{\delta,Q}(\zeta^L(0))$ for $r^L=r^{F_1}=r^{F_2}$ and $v^L_{\delta,Q'}(\zeta^L(0)) \geq v^L_{\delta,Q}(\zeta^L(0))$ for $-r^L=r^{F_1}=r^{F_2}$.
\end{theorem}

We remark that $Q^i$ is considered at least as informative as $Q'^i$ (see White and Harrington 1980, Zhang 2010). Theorem \ref{convex} implies that for a zero-sum game (a totally cooperative game), if (i) initially each follower has complete visibility of the other two agents' states and (ii) follower $i$'s value function is concave, then intuition regarding distorted information is still valid. Note that the concavity property is not required for $v^L$ in Theorem \ref{convex} and $v^{F_i}_{\delta,Q}$ is concave if follower $i$ uses an optimal response policy.
\section{Discussion}
\label{Implication}
\begin{figure}
	\begin{center}
		\includegraphics[width=0.75\textwidth]{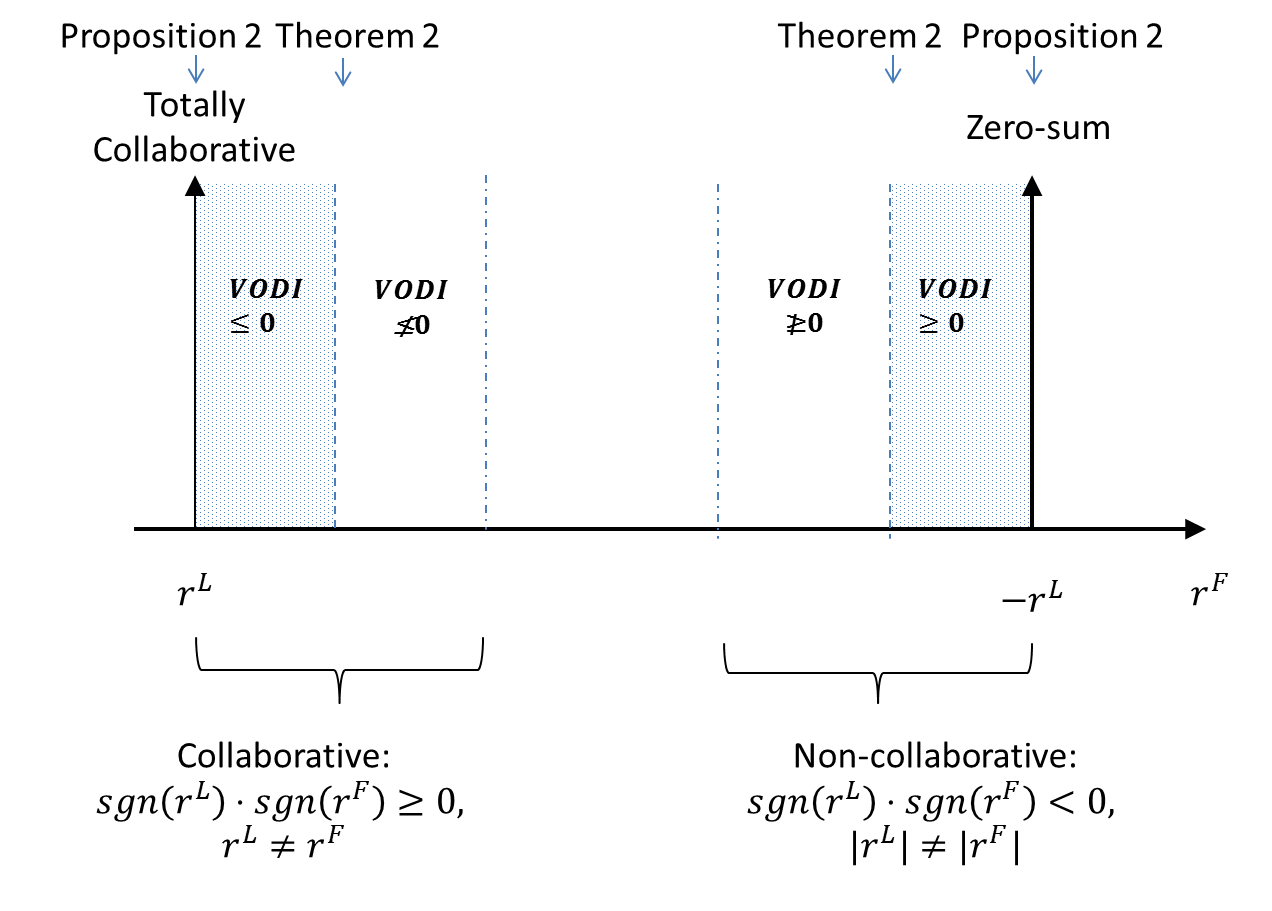}
		\caption{The value of distorted information (VODI) is negative (positive) for only a subset of collaborative (non-collaborative) games identified by Proposition \ref{1} and Theorem \ref{deviation}.}
		\label{explaination}
	\end{center}
\end{figure}

The analysis in the previous section shows that the intuitive understanding of distorted information could be invalid and the value of distorted information depends on the reward structure and the policies employed by all agents. Figure \ref{explaination} illustrates why the prevalent intuition may fail. Specifically, the value of distorted information is negative for the totally collaborative case and it is positive for the zero-sum case. Thus the sign of the value of distorted information will change as the game transitions from totally collaborative to zero-sum. However, the sign will not flip abruptly. During the transition, there are games whose value of distorted information is not positive (and not negative either) and the value of distorted information is negative (positive) for only a subset of collaborative (non-collaborative) games ``close" enough to the totally collaborative (zero-sum) case (Theorem \ref{deviation}). Similarly, manipulating distorted information can lead to a desired outcome only when the policies employed by agents are ``good" enough (see Figure \ref{randomPolicy}). These results suggest that decision makers (e.g., defense planners) should understand both the objective and the current policy (or practice) of each agent before modulating information, in order to obtain the desired benefit from intervention. \\

We remark that even if distorting (enhancing the accuracy of) information passing between adversaries (allies) can improve the leader's performance under the current reward structure and policies in use, this does not imply that the more investment in distorted information the better.  There is a threshold on investment, $c_1,c_2$ determined in Proposition \ref{Partition}, beyond which \textit{intelligent} followers (especially adversaries) may adjust their policies, leading to potentially unexpected degradation in the leader's performance. These insights are critical for making efficient targeted investments.

\subsection{When $\alpha^{1,0} \leq (\geq) \alpha^{0,1}$}
We remark that $\alpha^{0,1}, \alpha^{1,0}$ can be viewed as the marginal value to the leader of per unit distorted information (represented by $\epsilon_1, \epsilon_2$), assuming $\alpha^{0,0} + \alpha^{1,0}\epsilon_1 + \alpha^{0,1}\epsilon_2$ is a good approximation of $g^L$.  Below, we demonstrate that the marginal value of per unit distorted information is bounded by $\frac{2\beta M^L_\delta}{(1-\beta)^2}$, where $M^L_\delta = \max_{\zeta}|R^L_\delta(\zeta)|$.  We also present a set of conditions that guarantee $\alpha^{1,0} \leq (\geq) \alpha^{0,1}$. This result can be be useful in identifying efficient targeted investments.  We start with zero memory policies to illustrate the concepts without introducing additional and complicated notation. The challenges in analyzing general finite-memory policies are discussed in the Appendix.\\

Assume a zero-memory leader policy $\delta^L: S \rightarrow A^L$ and follower $i$'s zero-memory policy $\delta^{F_i}:S^F \times S^L_i \times Z^{F_i} \rightarrow A^{F_i}, i \neq j, i, j \in \{1,2\}$. Let $N=(|S||Z^{F_1}||Z^{F_2}|)$, and assume $\mu$ is a one-to-one, onto mapping from $\{\zeta\}$ to $\{1,2,...N\}$, where $\zeta = \{\zeta^L=\{s\}, \zeta^{F_1}=\{s^L_1,z^{F_1},s^F\}, \zeta^{F_2}=\{z^{F_2},s^L_2,s^F\} \} \equiv \{s,z^{F_1},z^{F_2}\}$. Thus, $\mu$ totally orders $\{\zeta\}$. A function $f: \{\zeta\} \rightarrow R$ is isotone in $\{s,s^L_2,s^L_1\}$ (with respect to $\mu$) if and only if $\mu(\{s,s^L_2,s^L_1\}) \leq \mu(\{s',s^{L'}_2,s^{L'}_1\})$ implies $f(\{s,s^L_2,s^L_1\}) \leq f(\{s',s^{L'}_2,s^{L'}_1\})$.\\

Define $q_\delta(k|\zeta) = \sum P(s'|s,\delta(\zeta)), o^{F_1}_\delta(k|\{s,s^L_2,s^L_1\}) =\sum (\sum_{z^{F_1}\neq s^{L}_2}\sigma_{s^{L}_2 ,z^{F_1}}P(s'|s,\delta(\{s,z^{F_1},s^L_1\}))),\\ o^{F_2}_\delta(k|\{s,s^L_2,s^L_1\}) = \sum( \sum_{z^{F_2}\neq s^{L}_1}\sigma_{s^{L}_1 ,z^{F_2}}P(s'|s,\delta(\{s,s^L_2,z^{F_2}\})))$, where the first sum of each term is over all $s'$ such that $\mu(\{s', s^{L'}_2, s^{L'}_1\}) \geq k$.
\begin{theorem}
	\label{quantify}
	$||\alpha^{0,1}|| \leq \frac{2\beta M^L_\delta}{(1-\beta)^2}$ and $||\alpha^{1,0}|| \leq \frac{2\beta M^L_\delta}{(1-\beta)^2}$. For $0 < \beta \leq 1/9$, $||\alpha^{k,l}|| \leq (\frac{8 \beta}{1-\beta})^{\max(k,l)}(\frac{M^L_\delta}{1-\beta})$, and for $1/9 < \beta < 1$, $||\alpha^{k,l}|| \leq (\frac{8\beta}{1-\beta})^{l+k}(\frac{M^L_\delta}{1-\beta})$. \\
	Furthermore, $\alpha^{1,0}\leq \alpha^{0,1}$ if the following conditions are satisfied:
	\begin{enumerate}[(i)]
		\item $R^L_\delta(\{s,s^L_2,s^L_1\})$ is isotone in $\{s,s^L_2,s^L_1\}$,
		\item $q_\delta(k|\zeta)$ is isotone in $\{s,s^L_2,s^L_1\}$,
		\item $\sum_{s^{L}_2 \neq z^{F_1} }\sigma_{s^{L}_2 ,z^{F_1}}R^L_\delta(\{s,z^{F_1},s^L_1\}) \leq \sum_{s^{L}_1 \neq z^{F_2} }\sigma_{s^{L}_1 ,z^{F_2}}R^L_\delta(\{s,s^L_2,z^{F_2}\}), \forall s$,
		\item $o^{F_1}_\delta(k|\{s,s^L_2,s^L_1\}) \leq o^{F_2}_\delta(k|\{s,s^L_2,s^L_1\}), \forall k, s$.
	\end{enumerate}
\end{theorem}

Generally speaking, if distorting information on channel 1 ($z^{F_1} \neq s^L_2$) tends to have greater impact on the leader's performance than distorting information on channel 2 ($z^{F_2} \neq s^L_1$), then it is likely that $\alpha^{1,0} \leq \alpha^{0,1}$. The $\alpha$ vectors for an example satisfying conditions in Theorem \ref{quantify} are listed in Table \ref{alphacompare}. The parameters for this example are included in the Appendix.
\begin{table}
	\begin{center}
		{\scriptsize
			\begin{tabular}{ | c | c | c | c | c | c | c | c | c | c | c | c |} \hline			

				$\zeta$ &$\alpha^{0,1}$ & $\alpha^{1,0}$ & $\zeta$
				&$\alpha^{0,1}$ & $\alpha^{1,0}$ & $\zeta$ &$\alpha^{0,1}$ & $\alpha^{1,0}$ & $\zeta$ &$\alpha^{0,1}$ & $\alpha^{1,0}$  \\
				\hline
				[0 0] [0 0] 0	&	-1.17	&	-3.85	&	[0 0] [0 1] 0	&	-0.62	&	-3.36	&	[1 0] [0 0] 0	&	0.37	&	-2.55	&	[1 0] [0 1] 0	&	0.11	&	-2.76	\\\hline
				[0 1] [0 0] 0	&	-1.33	&	-3.97	&	[0 1] [0 1] 0	&	-0.68	&	-3.4	&	[1 1] [0 0] 0	&	0.41	&	-2.44	&	[1 1] [0 1] 0	&	1.25	&	-1.82	\\\hline
				[0 0] [1 0] 0	&	-1.28	&	-3.93	&	[0 0] [1 1] 0	&	1.28	&	-1.36	&	[1 0] [1 0] 0	&	-0.18	&	-2.94	&	[1 0] [1 1] 0	&	0.86	&	-2.1	\\\hline
				[0 1] [1 0] 0	&	1.09	&	-1.57	&	[0 1] [1 1] 0	&	-0.4	&	-3.2	&	[1 1] [1 0] 0	&	0.53	&	-2.43	&	[1 1] [1 1] 0	&	0.31	&	-2.49	\\\hline
				[0 0] [0 0] 1	&	-0.92	&	-3.62	&	[0 0] [0 1] 1	&	-0.26	&	-3.05	&	[1 0] [0 0] 1	&	0.89	&	-2.12	&	[1 0] [0 1] 1	&	0.18	&	-2.73	\\\hline
				[0 1] [0 0] 1	&	-0.98	&	-3.66	&	[0 1] [0 1] 1	&	-0.43	&	-3.17	&	[1 1] [0 0] 1	&	0.31	&	-2.56	&	[1 1] [0 1] 1	&	1.51	&	-1.58	\\\hline
				[0 0] [1 0] 1	&	-0.87	&	-3.58	&	[0 0] [1 1] 1	&	1.37	&	-1.66	&	[1 0] [1 0] 1	&	0.06	&	-2.72	&	[1 0] [1 1] 1	&	1.36	&	-1.68	\\\hline
				[0 1] [1 0] 1	&	0.86	&	-2.18	&	[0 1] [1 1] 1	&	-0.04	&	-2.89	&	[1 1] [1 0] 1	&	0.72	&	-2.24	&	[1 1] [1 1] 1	&	0.55	&	-2.27	\\\hline		
		\end{tabular}}	
	\end{center}
	\caption{An example where $\alpha^{1,0} \leq \alpha^{0,1}$}
	\label{alphacompare}
\end{table}
Theorem \ref{quantify} in conjunction with Proposition \ref{Partition} provide managers with insights for making efficient targeted investments.
Assume that the joint policy $(\delta^L, \delta^{F_1}, \delta^{F_2})$ is such that $\alpha^{1,0} \leq \alpha^{0,1} \leq 0$, then it is optimal to modulate $\epsilon_1$ before $\epsilon_2$ for the higher marginal improvement, while $|\epsilon_1'-\epsilon_1|\leq c_1$. The leader can then continue to improve its performance by adjusting $\epsilon_2$ when $|\epsilon_1'-\epsilon_1|= c_1 $. Changing $\epsilon_1$ beyond this threshold should be approached with caution and be well justified.

\section{Conclusions and Future Research Directions}
\label{conclusion}
We have examined the value of distorted information in communication in a partially observable stochastic warfare game with three agents. Intuition suggests that the leader should benefit from less (more) accurate communication of information between adversarial (cooperative) followers. However, we demonstrated that such intuition is correct only under certain conditions. We investigated the range of deviations from these conditions under which the intuition remain intact. We further showed that when: (i) the reward structures are not perfectly aligned enough (or not exactly the opposite), or (ii) (even good) sub-optimal policies are employed, the communication of distorted information does not necessarily lead to the expected results.  We discussed why the prevalent intuition may not hold and proposed a holistic paradigm considering both the reward structures of agents and the policy employed by each agent for an improved understanding of the value of distorted information.  Without adopting such a holistic paradigm encompassing reward structures, policies, and information accuracy, decision makers could easily be led astray by their intuition.\\

The results of this paper provide meaningful guidance on the effective interventions and efficient investments.  We first identified under what conditions distorting information is beneficial to the leader (investigating ``if and when to invest").  We then obtained insights that can guide decision makers in making efficient targeted investments.  We bound the marginal improvement on the leader's performance per unit of distorted information, which can then be utilized to bound the marginal improvement per investment dollar (best ``bang per buck") if the relationship between the amount of investment and level of distorted information is known. We pointed out that more investment may not necessarily improve the leader's performance; the followers may adjust their policies when the communication quality is worse enough, which may not be favorable to the leader. \\

Although we presented our work in the context of warfare, our methodology and results may be applied to other critical domains such as cybersecurity, business competition, and political campaigns. For instance, our methodology could help an organization understand how to influence customers and discredit competitors in order to maximize margins and market share. Several researchers have identified the power of information in national politics and diplomatic relations.  It has been shown that spreading distorted information among a party's supporters may affect the popularity of a party and its policies, in the context of political contests or campaigns between two parties with opposing interests (Gul and Pesendorfer, 2012).  It would be be interesting to explore how our formulations and results could be applied to understand the role of information in these fields. As this research represents an initial investigation into the value of distorted information in communication to a leader, additional research on the impact of distorted information on online social networks would also be value adding. Indicatively, it would be interesting to explore the dependence of the value of distorted information on the various structures of communication networks to guide investment decisions.

\section{Acknowledgements}
The authors greatly benefited from email exchanges, discussions, and phone conversations with Maj. (ret) John W. Spencer, Chair of Urban Warfare Studies at the Modern War Institute and Co-Director of the Urban Warfare Project regarding the issues and challenges in the Battle of Mosul and wish to thank him for giving of his time so generously.

\section*{Appendix A: Technical Appendices}
\proof{PROOF OF PROPOSITION \ref{Partition}.}
The proof follows the similar line of reasoning as in Chang et al. (2015b) with $c_i = \frac{b_i(1-\beta)^2}{4\beta M^{F_i}}$, where $b_i = \min\{B(\mu^{F_i}): \mu^{F_i} \in \Pi^{F_i},  \mu^{F_i} \neq \delta^{F_i} \} > 0$, $B(\mu^{F_i}) = v^{F_i}_Q(\delta^L, \delta^{F_i}, \delta^{F_j})(\zeta^{F_i}) - v^{F_i}_Q(\delta^L, \mu^{F_i}, \delta^{F_j})(\zeta^{F_i})$ and $M^{F_i}=\max_s\max_a |r^{F_i}(s,a)|, i = 1,2$.
\endproof
\proof{PROOF OF LEMMA \ref{littleLemmaforConvergence}.}
(1) Clearly, $a^{0,0} \leq 1, a^{0,l} \leq c^l \leq (4c)^l, a^{k,0} \leq c^k \leq (4c)^k, \forall k, l \geq 0$. Assume $a^{k,l} \leq (4c)^{l+k}$ for $l+k \leq m$, then $a^{k+1,l} \leq c(a^{k,l}+a^{k+1,l-1})+2ca^{k,l-1} \leq (4c)^{k+l-1}(2c)(4c+1) \leq (4c)^{k+l+1}$ if $c \geq 1/4$. Similarly, $a^{k+1,l} \leq (4c)^{k+l+1}$.  \\
(2) W.L.O.G $k \leq l$. Initial conditions are clearly satisfied. Assume $a^{k,l} \leq (4c)^l$ for $l+k \leq m$. $a^{k+1,l} \leq c(a^{k,l}+a^{k+1,l-1})+2ca^{k,l-1} \leq c((4c)^l + (4c)^{\max(k+1,l-1)})+(2c)(4c)^{\max(k,l-1)}$. If $k \leq l-2$, $a^{k+1,l} \leq (4c)^{l-1}(c(4c+3)) \leq (4c)^{l} \leq (4c)^{\max(k+1,l)}$; if $k=l-1$, $a^{k+1,l} \leq (2c)(4c)^k(4c+1)\leq (4c)^{k+1} = (4c)^{\max(k+1,l)}$; if $k = l$, $a^{k+1,l} \leq (4c)^lc(3+4c) \leq (4c)^{l+1} =(4c)^{\max(k+1,l)}$ for $0 < c < 1/4$. Similarly, $a^{k,l+1} \leq c((4c)^{l+1}+(4c)^l)+(2c)(4c)^l \leq (4c)^{l+1} = (4c)^{\max(k,l+1)}$.	
\endproof

\proof{PROOF OF THEOREM \ref{thm1}.}
Step 1: Represent $v^L$ in terms of $g$.

Let $g^L_\delta$ be a solution to the equation
\begin{align*}
&g^L_\delta(\zeta(t,\tau)) = R^L_\delta(\zeta(t,\tau)) \\
&+ \beta \sum_{z(t+1)}\sum_{s(t+1)} P(z(t+1),s(t+1)|s(t),\delta(\zeta(t,\tau)))g^L_\delta(\bar{\zeta}(z(t+1),s(t+1), \delta(\zeta(t,\tau)), \zeta(t,\tau)) ),
\end{align*}
where $R^L_\delta(\zeta(t,\tau))= \sum_{a(t)} r^L(s(t), a(t))P(a(t)|\zeta(t,\tau))$. Then Proposition 2 in Chang et al. (2015a) and the fact that each agent uses a finite-memory policy show that $g^L_\delta$ is unique and $v^L_{\delta}(\zeta^L(t,\tau), y^L(t))= \sum_{\zeta^{F}(t,\tau)}P(\zeta^{F}(t,\tau)|\zeta^{L}(t))g^k_\delta(\zeta(t,\tau))$.\\
Step 2: Represent $g$ in terms of $\alpha$.\\
Let $g^L_0 = 0$ and $g^L_{n+1}(\zeta) = R^L_\delta(\zeta) + \beta \sum_{z'}\sum_{s'} P(z',s'|\zeta,\delta)g^L_n(\bar{\zeta}(z',s', \zeta,\delta(\zeta) )$, then $ \lim_{n\rightarrow \infty}||g^L_n - g^L|| = 0$. It follows from the definition of $Q^1$ and $Q^2$ that
\begin{align*}
&g^L_{n+1}(\zeta) = R^L_\delta(\zeta)  + \Lambda_n(\zeta)+\epsilon_1 \Delta_n(\zeta)
+\epsilon_2 \bar{\Delta}_n(\zeta)+\epsilon_1 \epsilon_2 \Theta_n(\zeta)
\end{align*}
where:
\begin{align*}
&\Lambda_n(\zeta)= \beta \sum_{s'}P(s'|s, \delta(\zeta))g^L_{n}( \bar{\zeta}[s',s^{L'}_2,s^{L'}_1,\zeta, \delta(\zeta)])\\
&\Delta_n(\zeta) = \beta \sum_{s'}P(s'|s, \delta(\zeta))\bigg[\sum_{s^{L'}_2 \neq z^{F_1'}}\sigma_{s^{L'}_2 ,z^{F_1'}}g^L_n(\bar{\zeta}[s',z^{F'_1},s^{L'}_1,\zeta, \delta(\zeta)] )-g^L_{n}(\bar{\zeta}[s',s^{L'}_2,s^{L'}_1,\zeta, \delta(\zeta)])
\bigg]\\
&\bar{\Delta}_n (\zeta) =	\beta \sum_{s'}P(s'|s, \delta(\zeta))\bigg[\sum_{s^{L'}_1 \neq z^{F_2'}}\sigma_{s^{L'}_1 ,z^{F_2'}}g^L_n(\bar{\zeta}[s',s^{L'}_2,z^{F'_2},\zeta, \delta(\zeta)])-g^L_{n}(\bar{\zeta}[s',s^{L'}_2,s^{L'}_1,\zeta, \delta(\zeta)])
\bigg]\\
&\Theta_n(\zeta) = \beta \sum_{s'}P(s'|s, \delta(\zeta))\bigg[\sum_{s^{L'}_1 \neq z^{F_2'}}\sum_{s^{L'}_2 \neq z^{F_1'}}\sigma_{s^{L'}_2 ,z^{F_1'}}\sigma_{s^{L'}_1 ,z^{F_2'}}  g^L_n(\bar{\zeta}[s',z^{F'_1},z^{F'_2},\zeta, \delta(\zeta)])\\
&-\sum_{s^{L'}_2 \neq z^{F_1'}}\sigma_{s^{L'}_2 ,z^{F_1'}}g^L_n(\bar{\zeta}[s',z^{F'_1},s^{L'}_1,\zeta, \delta(\zeta)] )-\sum_{s^{L'}_1 \neq z^{F_2'}}\sigma_{s^{L'}_1 ,z^{F_2'}}g^L_n(\bar{\zeta}[s',s^{L'}_2,z^{F'_2},\zeta, \delta(\zeta)])\\
&+g^L_{n}(\bar{\zeta}[s',s^{L'}_2,s^{L'}_1,\zeta, \delta(\zeta)])
\bigg].
\end{align*}

Then it is straightforward to show that $g_n^L(\zeta) = \sum_{k=0}^{n-1} \sum_{l=0}^{n-1} \epsilon_1^k \epsilon_2^l \alpha_n^{k,l}(\zeta) $, where:
\begin{align*}
&\alpha^{0,0}_{n+1}  = R^L_\delta + \Lambda_n^{0,0},\hspace{5pt}\alpha^{0,l}_{n+1} = \bar{\Delta}^{0,l-1}_n + \Lambda_n^{0,l},\hspace{5pt} \alpha^{k,0}_{n+1}  = \Delta^{k-1,0}_n   + \Lambda_n^{k,0} ,\hspace{5pt}k,l = 1,...,n-1;  \\	
&\alpha^{0,n}_{n+1}  = \bar{\Delta}^{0,n-1}_n,\hspace{5pt}\alpha^{n,0}_{n+1} = \Delta^{n-1,0}_n; \\
&\alpha^{k,n}_{n+1} = \bar{\Delta}^{k,n-1}_n+ \Theta_n^{k-1,n-1},\hspace{5pt}\alpha^{n,l}_{n+1} = \Delta^{n-1,l}_n + \Theta_n^{n-1,l-1},\hspace{5pt}k,l=1,...,n-1; \\
&\alpha^{k,l}_{n+1} = \Delta^{k-1,l}_n + \bar{\Delta}^{k,l-1}_n + \Theta^{k-1,l-1}_n+ \Lambda_n^{k,l}, \hspace{5pt}k,l = 1,...,n-1;\\
&\alpha^{n,n}_{n+1} =  \Theta^{n-1,n-1}_n.\\
\mbox{and}\\
&\Lambda_n^{k,l}(\zeta) = \beta\sum_{s'}P(s'|s,\delta(\zeta))\alpha^{k,l}_n(\bar{\zeta}[s',s^{L'}_2,s^{L'}_1,\zeta, \delta(\zeta)]);\\
&\Delta^{k,l}_n(\zeta) = -\Lambda_n^{k,l}(\zeta) + \beta \sum_{s'}P(s'|s, \delta(\zeta))\sum_{s^{L'}_2 \neq z^{F_1'} } \sigma_{s^{L'}_2 ,z^{F_1'}} \alpha^{k,l}_n(\bar{\zeta}[s',z^{F'_1},s^{L'}_1,\zeta, \delta(\zeta)] ) , k, l \geq 0;\\
&\bar{\Delta}^{k,l}_n(\zeta) = -\Lambda_n^{k,l}(\zeta) + \beta \sum_{s'}P(s'|s, \delta(\zeta))\sum_{s^{L'}_1 \neq z^{F_2'} } \sigma_{s^{L'}_1 ,z^{F_2'}} \alpha^{k,l}_n(\bar{\zeta}[s',s^{L'}_2,z^{F'_2},\zeta, \delta(\zeta)]) ,k, l \geq 0;\\
&\Theta^{k,l}_n(\zeta) =\Lambda_n^{k,l}(\zeta)+ \beta \sum_{s'}P(s'|s,\delta(\zeta))\bigg[\sum_{z^{F_1'} \neq s^{L'}_2}\sum_{z^{F_2'} \neq s^{L'}_1} \sigma_{s^{L'}_2 ,z^{F_1'}}\sigma_{s^{L'}_1 ,z^{F_2'}}\alpha^{k,l}_n(\bar{\zeta}[s',z^{F'_1},z^{F'_2},\zeta, \delta(\zeta)])\\
& -\sum_{z^{F_1'} \neq s^{L'}_2} \sigma_{s^{L'}_2 ,z^{F_1'}}\alpha^{k,l}_n(\bar{\zeta}[s',z^{F'_1},s^{L'}_1,\zeta, \delta(\zeta)] )
-\sum_{z^{F_2'} \neq s^{L'}_1} \sigma_{s^{L'}_1 ,z^{F_2'}}\alpha^{k,l}_n(\bar{\zeta}[s',s^{L'}_2,z^{F'_2},\zeta, \delta(\zeta)]) \bigg], k,l \geq 1.
\end{align*}
Assuming $\lim_{n\rightarrow \infty} \sum_{k=0}^{n-1} \sum_{l=0}^{n-1} \epsilon_1^k \epsilon_2^l \alpha_n^{k,l}(\zeta)$ exists (the absolute convergence is proved in step 3), taking the limit gives: 	\begin{align*}
&\alpha^{0,0} = R^L_\delta +\Lambda^{0,0},\hspace{5pt}\alpha^{0,l} = \bar{\Delta}^{0,l-1}+\Lambda^{0,l},\hspace{5pt}\alpha^{k,0}= \Delta^{k-1,0}+ \Lambda^{k,0}, \\
&\alpha^{k,l} = \Delta^{k-1,l} + \bar{\Delta}^{k,l-1} + \Theta^{k-1,l-1}+ \Lambda^{k,l},\hspace{5pt}k,l = 1,...,n-1. \\
&\mbox{and}\\
&\Lambda^{k,l}(\zeta) = \beta\sum_{s'}P(s'|s,\delta(\zeta))\alpha^{k,l}(\bar{\zeta}[s',s^{L'}_2,s^{L'}_1,\zeta, \delta(\zeta)]), k,l \geq 0;\\
&\Delta^{k,l}(\zeta) =-\Lambda^{k,l}(\zeta)+ \beta \sum_{s'}P(s'|s, \delta(\zeta)) \sum_{s^{L'}_2 \neq z^{F_1'} } \sigma_{s^{L'}_2 ,z^{F_1'}} \alpha^{k,l}(\bar{\zeta}[s',z^{F'_1},s^{L'}_1,\zeta, \delta(\zeta)] ) , k, l \geq 0;\\
&\bar{\Delta}^{k,l}(\zeta) =-\Lambda^{k,l}(\zeta)+ \beta \sum_{s'}P(s'|s, \delta(\zeta)) \sum_{s^{L'}_1 \neq z^{F_2'} } \sigma_{s^{L'}_1 ,z^{F_2'}} \alpha^{k,l}(\bar{\zeta}[s',s^{L'}_2,z^{F'_2},\zeta, \delta(\zeta)])  ,k, l \geq 0;\\
&\Theta^{k,l}(\zeta) = \Lambda^{k,l}(\zeta)+\beta \sum_{s'}P(s'|s,\delta(\zeta))\bigg[\sum_{z^{F_1'} \neq s^{L'}_2}\sum_{z^{F_2'} \neq s^{L'}_1} \sigma_{s^{L'}_2 ,z^{F_1'}}\sigma_{s^{L'}_1 ,z^{F_2'}}\alpha^{k,l}(\bar{\zeta}[s',z^{F'_1},z^{F'_2},\zeta, \delta(\zeta)])\\
& -\sum_{z^{F_1'} \neq s^{L'}_2} \sigma_{s^{L'}_2 ,z^{F_1'}}\alpha^{k,l}(\bar{\zeta}[s',z^{F'_1},s^{L'}_1,\zeta, \delta(\zeta)] )
-\sum_{z^{F_2'} \neq s^{L'}_1} \sigma_{s^{L'}_1 ,z^{F_2'}}\alpha^{k,l}(\bar{\zeta}[s',s^{L'}_2,z^{F'_2},\zeta, \delta(\zeta)]) \bigg], k,l \geq 1.
\end{align*}
Step 3: Convergence.\\
Define $M^L_\delta = \max_{\zeta}|R^L_\delta(\zeta)|$ and let $||u||=\max_\zeta|u(\zeta)|$ for a vector $u$. We note that $||\alpha^{0,0}|| \leq \frac{M^L_\delta}{1-\beta}$, $||\Delta^{k,l}|| \leq 2\beta ||\alpha^{k,l}||$, $||\bar{\Delta}^{k,l}|| \leq 2\beta ||\alpha^{k,l}||$, $||\alpha^{0,l}|| \leq \beta||\alpha^{0,l}|| + ||\bar{\Delta}^{0,l-1}||$, and $||\alpha^{k,0}|| \leq \beta||\alpha^{k,0}|| + ||\Delta^{k-1,0}||$. Hence, $||\alpha^{0,l}|| \leq (\frac{2\beta}{1-\beta})^l (\frac{M^L_\delta}{1-\beta})$ and $||\alpha^{k,0}|| \leq (\frac{2\beta}{1-\beta})^k (\frac{M^L_\delta}{1-\beta})$. Furthermore, $||\Theta^{k,l}|| \leq 4\beta ||\alpha^{k,l}|| $ and $||\alpha^{k,l}|| \leq \beta||\alpha^{k,l}|| + ||\Delta^{k-1,l}|| + ||\bar{\Delta}^{k,l-1}||+ ||\Theta^{k-1, l-1}||$ imply $||\alpha^{k,l}|| \leq \frac{2\beta}{1-\beta}(||\alpha^{k-1,l}|| + ||\alpha^{k,l-1}||) + \frac{4\beta}{1-\beta}||\alpha^{k-1,l-1}||$. Thus, $||\alpha^{k,l}|| \leq (\frac{8\beta}{1-\beta})^{l+k}(\frac{M^L_\delta}{1-\beta})$ for $1/9 \leq \beta <1 $ and $||\alpha^{k,l}|| \leq (\frac{8\beta}{1-\beta})^{\max(k,l)}(\frac{M^L_\delta}{1-\beta})$ for $0 < \beta \leq 1/9 $ by Lemma \ref{littleLemmaforConvergence}. For $1/9 \leq \beta < 1$, $||\sum_{k=0}^\infty \sum_{l=0}^\infty \epsilon_1^k \epsilon_2^l \alpha^{k,l}|| \leq \frac{M^L_\delta}{1-\beta}\sum_{k=0}^\infty \sum_{l=0}^\infty \epsilon_1^k \epsilon_2^l  (\frac{8\beta}{1-\beta})^{l+k} < \infty$ if $0<\epsilon_1, \epsilon_2 < \frac{1-\beta}{8\beta}$; for $0 < \beta < 1/9, ||\sum_{k=0}^\infty \sum_{l=0}^\infty \epsilon_1^k \epsilon_2^l \alpha^{k,l}|| \leq \frac{M^L_\delta}{1-\beta}\sum_{k=0}^\infty \sum_{l=0}^\infty \epsilon_1^k \epsilon_2^l  (\frac{8\beta}{1-\beta})^{\max\{k,l\}} \leq \frac{M^L_\delta}{1-\beta}\sum_{k=0}^\infty \sum_{l=0}^\infty \epsilon_1^k \epsilon_2^l <\infty$ for $\forall 0 < \epsilon_1, \epsilon_2 < 1$.
\endproof
\proof{PROOF OF PROPOSITION \ref{1}.}
Let $\{v_n^F\}$ and $\{\delta^{F*}_n\}$ be such that
$$v_{n+1}^F(s) = \max_{a^F} \bigg\{ R^F(s,\delta^{L}(s), a^F) + \beta \sum_{s'}P(s'|s,\delta^L(s),a^F) v_n^F(s') \bigg\} $$
where $v^F_0 = 0$. Then, $\lim_{n \rightarrow \infty} ||v^{F*} - v^F_n|| = 0$. If $r^L = r^F$, it follows by induction and the definition of $\alpha^{0,0}$ that $\alpha^{0,0}_{n+1}(\bar{\zeta}[s, \zeta^{L}, \delta^{L}(\zeta^{L})], \bar{\zeta}[s, \zeta^{F_1}, \delta^{F_1*}(\zeta^{F_1})], \bar{\zeta}[s, \zeta^{F_2}, \delta^{F_2*}(\zeta^{F_2})])=R^F(s, \delta^L(s) , \delta^{F*}_n(s)) + \beta \sum_{s'}P(s'|s,\delta^L(s), \delta^{F*}_n(s))v^F_n(s')= v^{F}_{n+1}(s)$.  Hence, $$\alpha^{0,0}(\bar{\zeta}[s, \zeta^{L}, \delta^{L}(\zeta^{L})], \bar{\zeta}[s, \zeta^{F_1}, \delta^{F_1*}(\zeta^{F_1})], \bar{\zeta}[s, \zeta^{F_2}, \delta^{F_2*}(\zeta^{F_2})])=v^{F*}(s).$$ Now if $z^{F_i} \neq s^L_j$, $\alpha^{0,0}(\bar{\zeta}[s, \zeta^{L}, \delta^{L}(\zeta^{L})], \bar{\zeta}[z^{F_i},s^F,s^L_j, \zeta^{F_i}, \delta^{F_i*}(\zeta^{F_i})], \bar{\zeta}[s, \zeta^{F_j}, \delta^{F_j*}(\zeta^{F_j})]) = \\R^F(s,\delta^L(s), \delta^{F_i*}(\{z^{F_i},s^F,s^L_j\}), \delta^{F_j*}(s))+ \beta \sum_{s'}P(s'|s,\delta^L(s), \delta^{F_i*}(\{z^{F_i},s^F,s^L_j\}), \delta^{F_j*}(s))v^{F*}(s') \leq R^F(s,\delta^L(s), \delta^{F_i*}(s), \delta^{F_j*}(s))+ \beta \sum_{s'}P(s'|s,\delta^L(s), \delta^{F_i*}(s), \delta^{F_j*}(s))v^{F*}(s')= v^{F*}(s)$. Hence, $\Delta^{0,0} \leq 0$, $\bar{\Delta}^{0,0} \leq 0$ and $\alpha^{1,0} \leq 0, \alpha^{0,1} \leq 0$ from the definition of $\alpha^{1,0}, \alpha^{1,0}, \Delta^{0,0}, \bar{\Delta}^{0,0}$ in terms of $\alpha^{0,0}$, and $v^L_{Q'}(\zeta^L, y^L) \leq v^L_Q(\zeta^L, y^L)$ from Theorem \ref{thm1}. The proof for the case where $r^L = -r^F$ is similar.
\endproof
\proof{PROOF OF LEMMA \ref{distanceProposition}.}
(a) Clearly, $d(\delta,\delta') \leq 2$. Let $U' = \{u: u \in U, u(\bar{\zeta}(s,s^L_2,s^L_1,\zeta,\delta(\zeta)))= u(\bar{\zeta}(s,s^L_2,s^L_1,\zeta,\delta'(\zeta))),  \forall \zeta, \forall \delta, \delta' \in \Pi \}$. Thus, $U' \subset U$.
\begin{align*}
&d(\delta,\delta')\\
& \geq
\sup_{||u|| \leq 1, u \in U'}||\sum_{s'}P(s'|s,\delta(\zeta))u(\bar{\zeta}(s',s^{L'}_2,s^{L'}_1,\zeta,\delta(\zeta)))
-\sum_{s'}P(s'|s,\delta'(\zeta))u(\bar{\zeta}(s',s^{L'}_2,s^{L'}_1,\zeta,\delta'(\zeta)))||\\
&=	\sup_{||u|| \leq 1, u \in U'}||\sum_{s'}P(s'|s,\delta(\zeta))u(\bar{\zeta}(s',s^{L'}_2,s^{L'}_1,\zeta,\delta(\zeta)))
-\sum_{s'}P(s'|s,\delta'(\zeta))u(\bar{\zeta}(s',s^{L'}_2,s^{L'}_1,\zeta,\delta(\zeta)))||\\
&=\sup_{||u|| \leq 1, u \in U'}||\sum_{s'}[P(s'|s,\delta(\zeta))-P(s'|s,\delta'(\zeta))]u(\bar{\zeta}(s',s^{L'}_2,s^{L'}_1,\zeta,\delta(\zeta)))||\\
&=||P(\delta) - P(\delta')||.
\end{align*}

(b) $U$ is the set of all bounded, real-valued functions generated by all possible zero memory policy tuples. Hence, $\forall u \in U, u(\zeta(s,s^{L}_2,s^{L}_1,\zeta,\delta(\zeta))) = u(\zeta(s,s^{L}_2,s^{L}_1,\zeta',\delta'(\zeta'))),  \forall \zeta\neq \zeta', \delta\neq \delta'$. The result follows from the proof in (a).
\endproof
\begin{lemma}
	\label{distance}
	\begin{enumerate}[(a)]
		\item[]
		\item $||\sum_{s'}P(s'|s,\delta(\zeta))u(\bar{\zeta}(s',s^{L'}_2,s^{L'}_1,\zeta,\delta(\zeta)))
		-\sum_{s'}P(s'|s,\delta'(\zeta))u(\bar{\zeta}(s',s^{L'}_2,s^{L'}_1,\zeta,\delta'(\zeta)))|| \leq d(\delta,\delta')||u||$,
		\item $||\sum_{s'}P(s'|s,\delta(\zeta))u(\bar{\zeta}(s',s^{L'}_2,s^{L'}_1,\zeta,\delta(\zeta)))
		-\sum_{s'}P(s'|s,\delta'(\zeta))w(\bar{\zeta}(s',s^{L'}_2,s^{L'}_1,\zeta,\delta'(\zeta)))|| \leq d(\delta,\delta')||u||+ ||u(\bar{\zeta}(s',s^{L'}_2,s^{L'}_1,\zeta,\delta'(\zeta)))-w(\bar{\zeta}(s',s^{L'}_2,s^{L'}_1,\zeta,\delta'(\zeta))) || $.
	\end{enumerate}	
\end{lemma}
\proof{PROOF OF LEMMA \ref{distance}.}
(a) is from the definitions of $d(\delta,\delta')$ and $P$ is a linear operator.\\
(b)
\begin{align*}
&||\sum_{s'}P(s'|s,\delta(\zeta))u(\bar{\zeta}(s',s^{L'}_2,s^{L'}_1,\zeta,\delta(\zeta)))
-\sum_{s'}P(s'|s,\delta'(\zeta))w(\bar{\zeta}(s',s^{L'}_2,s^{L'}_1,\zeta,\delta'(\zeta)))|| \\
&=||\sum_{s'}P(s'|s,\delta(\zeta))u(\bar{\zeta}(s',s^{L'}_2,s^{L'}_1,\zeta,\delta(\zeta)))
-\sum_{s'}P(s'|s,\delta'(\zeta))u(\bar{\zeta}(s',s^{L'}_2,s^{L'}_1,\zeta,\delta'(\zeta))) \\
&+ \sum_{s'}P(s'|s,\delta'(\zeta))u(\bar{\zeta}(s',s^{L'}_2,s^{L'}_1,\zeta,\delta'(\zeta))) -\sum_{s'}P(s'|s,\delta'(\zeta))w(\bar{\zeta}(s',s^{L'}_2,s^{L'}_1,\zeta,\delta'(\zeta))) ||\\
&\leq d(\delta,\delta')||u|| + ||u(\bar{\zeta}(s',s^{L'}_2,s^{L'}_1,\zeta,\delta'(\zeta)))-w(\bar{\zeta}(s',s^{L'}_2,s^{L'}_1,\zeta,\delta'(\zeta))) ||.
\end{align*}
\endproof
\proof{PROOF OF THEOREM \ref{deviation}.}
Let \\$u(\bar{\zeta}(s',s^{L'}_2,s^{L'}_1, \zeta, \rho(\zeta))) = \sum_{z^{F_1'}\neq s^{L'}_2}\sigma_{s^{L'}_2,z^{F_1'}}\alpha^{0,0}_{\delta^*}(\bar{\zeta}(s', z^{F_1'},s^{L'}_1, \zeta, \rho(\zeta))) - \alpha^{0,0}_{\delta^*}(\bar{\zeta}(s',s^{L'}_2,s^{L'}_1, \zeta, \rho(\zeta)))$,\\
and \\$w(\bar{\zeta}(s',s^{L'}_2,s^{L'}_1, \zeta, \rho(\zeta))) = \sum_{z^{F_1'}\neq s^{L'}_2}\sigma_{s^{L'}_2,z^{F_1'}}\alpha^{0,0}_\delta(\bar{\zeta}(s', z^{F_1'},s^{L'}_1, \zeta, \rho(\zeta))) - \alpha^{0,0}_\delta(\bar{\zeta}(s',s^{L'}_2,s^{L'}_1, \zeta, \rho(\zeta)))$.
\begin{align*}
&||\Delta_\delta^{0,0} - \Delta_{\delta^*}^{0,0}||\\
&=\beta|| \sum_{s'}P(s'|s,\delta(\zeta))[\sum_{z^{F_1'}\neq s^{L'}_2}\sigma_{s^{L'}_2,z^{F_1'}}\alpha^{0,0}_\delta(\bar{\zeta}(s', z^{F_1'},s^{L'}_1, \zeta, \delta(\zeta))) - \alpha^{0,0}_\delta(\bar{\zeta}(s',s^{L'}_2,s^{L'}_1, \zeta, \delta(\zeta)))]\\
&- \sum_{s'}P(s'|s,\delta^*(\zeta))[\sum_{z^{F_1'}\neq s^{L'}_2}\sigma_{s^{L'}_2,z^{F_1'}}\alpha^{0,0}_{\delta^*}(\bar{\zeta}(s', z^{F_1'},s^{L'}_1, \zeta, \delta^*(\zeta))) - \alpha^{0,0}_{\delta^*}(\bar{\zeta}(s',s^{L'}_2,s^{L'}_1, \zeta, \delta^*(\zeta)))]||\\
&= \beta||\sum_{s'}P(s'|s,\delta^*(\zeta))u(\bar{\zeta}(s',s^{L'}_2,s^{L'}_1, \zeta, \delta^*(\zeta)))- \sum_{s'}P(s'|s,\delta(\zeta))w(\bar{\zeta}(s',s^{L'}_2,s^{L'}_1, \zeta, \delta(\zeta))) ||\\
& \leq \beta d(\delta,\delta^*) ||u(\bar{\zeta}(s',s^{L'}_2,s^{L'}_1, \zeta, \delta^*(\zeta)))|| 	+ \beta\left (||\alpha^{0,0}_{\delta^*}(\bar{\zeta}(s',s^{L'}_2,s^{L'}_1, \zeta, \delta(\zeta)))- \alpha^{0,0}_{\delta}(\bar{\zeta}(s',s^{L'}_2,s^{L'}_1, \zeta, \delta(\zeta)))||\right ) \\
& + \beta\left (||\sum_{z^{F_1'}\neq s^{L'}_2}\sigma_{s^{L'}_2,z^{F_1'}}\alpha^{0,0}_{\delta^*}(\bar{\zeta}(s', z^{F_1'},s^{L'}_1, \zeta, \delta(\zeta))) - \sum_{z^{F_1'}\neq s^{L'}_2}\sigma_{s^{L'}_2,z^{F_1'}}\alpha^{0,0}_\delta(\bar{\zeta}(s', z^{F_1'},s^{L'}_1, \zeta, \delta(\zeta)))||\right ) \\
&\mbox{where the last inequality comes from Lemma \ref{distance}.}\\
&||\alpha_\delta^{1,0} - \alpha_{\delta^*}^{1,0}|| \leq ||\Delta_\delta^{0,0} - \Delta_{\delta^*}^{0,0}|| \\& +
\beta ||\sum_{s'}P(s'|s,\delta^*(\zeta))\alpha^{1,0}_{\delta^*}(\bar{\zeta}(s',s^{L'}_2,s^{L'}_1, \zeta, \delta^*(\zeta)))-\sum_{s'}P(s'|s,\delta(\zeta))\alpha^{1,0}_{\delta}(\bar{\zeta}(s',s^{L'}_2,s^{L'}_1, \zeta, \delta(\zeta)))|| \\
&\leq ||\Delta_\delta^{0,0} - \Delta_{\delta^*}^{0,0}|| + \beta\left(d (\delta,\delta^*)||\alpha^{1,0}_{\delta^*}|| + ||\alpha_\delta^{1,0}(\bar{\zeta}(s',s^{L'}_2,s^{L'}_1, \zeta, \delta(\zeta))) - \alpha_{\delta^*}^{1,0}(\bar{\zeta}(s',s^{L'}_2,s^{L'}_1, \zeta, \delta(\zeta)))||\right)\\
&\leq ||\Delta_\delta^{0,0} - \Delta_{\delta^*}^{0,0}|| + \beta\left(d (\delta,\delta^*)||\alpha^{1,0}_{\delta^*}|| + ||\alpha_\delta^{1,0}- \alpha_{\delta^*}^{1,0}||\right).\\
&||\alpha_\delta^{1,0}- \alpha_{\delta^*}^{1,0}||\leq \frac{1}{1-\beta} ||\Delta_\delta^{0,0} - \Delta_{\delta^*}^{0,0}|| + \frac{\beta}{1-\beta}\left(d (\delta,\delta^*)||\alpha^{1,0}_{\delta^*}|| \right)\\
&\leq \frac{\beta}{1-\beta} d(\delta,\delta^*) \eta^1 + \frac{\beta}{1-\beta}\left (||\alpha^{0,0}_{\delta^*}(\bar{\zeta}(s',s^{L'}_2,s^{L'}_1, \zeta, \delta(\zeta)))- \alpha^{0,0}_{\delta}(\bar{\zeta}(s',s^{L'}_2,s^{L'}_1, \zeta, \delta(\zeta)))||\right ) \\
& + \frac{\beta}{1-\beta}\left (||\sum_{z^{F_1'}\neq s^{L'}_2}\sigma_{s^{L'}_2,z^{F_1'}}\left[\alpha^{0,0}_{\delta^*}(\bar{\zeta}(s', z^{F_1'},s^{L'}_1, \zeta, \delta(\zeta))) - \alpha^{0,0}_\delta(\bar{\zeta}(s', z^{F_1'},s^{L'}_1, \zeta, \delta(\zeta)))\right]||\right )
\end{align*}
The result follows by the fact that if $||\alpha_\delta^{1,0} - \alpha_{\delta^*}^{1,0}|| < h_1$, then $\alpha_\delta^{1,0}  < 0$.
The proof for $\alpha^{0,1}$ is similar.
\endproof
\proof{PROOF OF COROLLARY \ref{deviationcorollary}.}
Under zero-memory policies,
\begin{align*}
&||\alpha_\delta^{1,0} - \alpha_{\delta^*}^{1,0}||  \\&\leq \frac{\beta}{1-\beta}d(\delta,\delta^*)\eta^1 +\frac{\beta}{1-\beta}\left (||\sum_{z^{F_1}\neq s^{L}_2}\sigma_{s^{L}_2,z^{F_1}}\left[\alpha^{0,0}_{\delta^*}(\{s, z^{F_1},s^L_1\}) - \alpha^{0,0}_\delta(\{s, z^{F_1},s^L_1\})\right]||+||\alpha^{0,0}_{\delta^*}- \alpha^{0,0}_{\delta}||\right ) \\
&\leq \frac{\beta}{1-\beta}d(\delta,\delta^*)\eta^1 +\frac{2\beta}{1-\beta}||\alpha^{0,0}_{\delta^*} - \alpha^{0,0}_\delta||.
\end{align*}
By Lemma \ref{distanceProposition}(b), we have $\alpha^{1,0} \leq 0$
if $\eta^1 ||P(\delta)-P(\delta^*)||+2||\alpha^{0,0}_{\delta^*} - \alpha^{0,0}_\delta|| \leq \frac{1-\beta}{\beta}h_1$. The proof for $\alpha^{0,1}$ is similar.
\endproof

\proof{PROOF OF PROPOSITION \ref{2}.}
Let $a^F = a^*, \forall \zeta^F$ and define an operator $H$ on vector $v$ $$[Hv](\zeta^L)=R^L(\zeta^L,\delta^L(\zeta^L),a^*) + \beta \sum_{s'}P(s'|s, \delta^L(\zeta^L),a^*)v(\bar{\zeta}[s',\zeta^L, \delta^{L}(\zeta^{L})])$$
Define the sequence $\{\tilde{\alpha}^{0,0}_n \}$ as $\tilde{\alpha}^{0,0}_{n+1} = H \tilde{\alpha}^{0,0}_n$, where $\tilde{\alpha}^{0,0}_0 = 0$. Thus, $\tilde{\alpha}^{0,0}= \lim_{n\rightarrow \infty}\tilde{\alpha}^{0,0}_n$ exists, and $\tilde{\alpha}^{0,0}(\zeta^{L})=[H\tilde{\alpha}^{0,0}](\zeta^{L})$. Assume $\alpha^{0,0}_n(\zeta^L, \zeta^F) = \tilde{\alpha}^{0,0}_n(\zeta^L), \forall \zeta^L, \forall \zeta^F$, then from the definition of $\alpha^{0,0}_{n+1}(\zeta^L,\zeta^F) $ in the proof of Theorem \ref{thm1} and by induction, $\alpha^{0,0}(\zeta^L,\zeta^F)= \lim_{n\rightarrow \infty} \alpha_{n}^{0,0}(\zeta^L,\zeta^F) =  \tilde{\alpha}^{0,0}(\zeta^L), \forall \zeta^L, \forall \zeta^F$. Hence, $\Delta^{0,0}=0$ and $\bar{\Delta}^{0,0}=0$ by the definitions of $\Delta^{0,0}$ and $\bar{\Delta}^{0,0}$. Accordingly, $\alpha^{0,1} = \alpha^{1,0} = 0$, $\Delta^{0,1}=\bar{\Delta}^{1,0}=0$, and $\Theta^{0,0} = 0$, which further imply $\alpha^{1,1} = 0$. Now assume $\alpha^{k,0} = \alpha^{0,l} =0, \forall k, l = m$, then $\Delta^{k,0}=0$, $\bar{\Delta}^{0,l}=0$ and as a result $\alpha^{k+1,0} = \alpha^{0,l+1} = 0$. Assume $\forall k+l \leq m,\alpha^{k,l} = 0$, then $\Delta^{k,l} = \bar{\Delta}^{k+1,l-1} = 0, \Theta^{k,l-1} = \Theta^{k-1,l} = 0$. Hence, $\alpha^{k+1,l} = \Delta^{k,l} + \bar{\Delta}^{k+1,l-1}+\Theta^{k,l-1}+\Lambda^{k+1,l} = \Lambda^{k+1,l}$. Thus, $\alpha^{k+1,l} = 0$. Similarly, $\alpha^{k,l+1} = 0$. Hence, $v^L(\zeta^L, y^L) = \sum_{\zeta^F} y^L(\zeta^F) \alpha^{0,0}(\zeta^L,\zeta^F)$ and the leader's value function is independent of $\epsilon_1$ and $\epsilon_2$.
\endproof
\proof{PROOF OF THEOREM \ref{convex}.}
It is easy to show that
\begin{align*}
&v^{F_i}_{\delta}(\zeta^{F_i}, y^{F_i}) =\sum_{\zeta^{F_j}, \zeta^L}R^{F_i}_\delta(\zeta)y(\zeta^{F_j}, \zeta^L) + \beta \sum_{z^{'}}\sum_{s^{'}}\sum_{\zeta^{F_j}}\sum_{\zeta^L}P(z',s'|s,\delta(\zeta))y(\zeta^{F_j}, \zeta^L) \notag\\
&\times v^{F_i}_{\delta}(\bar{\zeta}(z^{F'_i},s^{L'}_i,s^{F'},\delta^{F_i}(\zeta^{F_i}), \zeta^{F_i}), \lambda(z^{F'}_i,s^{L'}_i,s^{F'},\delta^{F_i}(\zeta^{F_i}), \zeta^{F_i}, y^{F_i})),
\end{align*}
where
$\lambda(z^{F'}_i,s^{L'}_i,s^{F'},\delta^{F_i}(\zeta^{F_i}), \zeta^{F_i}, y^{F_i})$ is the stochastic array with scalar element \\ $P(\zeta^{F_j'}, \zeta^{L'}| z^{F'_i},s^{L'}_i,s^{F'},\delta^{F_i}(\zeta^{F_i}), \zeta^{F_i}) = \frac{P(z^{F'_i},s^{F'},s^{L'}_i, \zeta^{F'_j},\zeta^{L'}|\zeta^{F_i}, y^{F_i})}{\sum_{z^{F_j'}}\sum_{s^{L'}_j}\sum_{\zeta^L}\sum_{\zeta^{F_j}}P(z',s'|s,\delta(\zeta))y(\zeta^{F_j},\zeta^L)}$, \\assuming $\sum_{z^{F_j'}}\sum_{s^{L'}_j}\sum_{\zeta^L}\sum_{\zeta^{F_j}}P(z',s'|s,\delta(\zeta))y(\zeta^{F_j},\zeta^L) \neq 0$.

Thus $v^{F_i}_{\delta,Q'}(\zeta^{F_i}(t,\tau),y^{F_i}(t)) \leq v^{F_i}_{\delta,Q}(\zeta^{F_i}(t,\tau),y^{F_i}(t))$ follows directly from White and Harrington (1980). $P(\zeta^{L,F_j}(0)|\zeta^{F_i}(0))=1$ implies $v^{F_i}_{\delta,Q}(\zeta^{F_i}(0))=g^{F_i}_Q(\zeta(0))$. The results follow the fact that $g^L_Q=g^{F_1}_Q=g^{F_2}_Q$ for $r^L = r^{F_1}=r^{F_2}$ and $-g^L_Q=g^{F_1}_Q=g^{F_2}_Q$ for $-r^L = r^{F_1}=r^{F_2}$.
\endproof
\proof{PROOF OF THEOREM \ref{quantify}.}
The bounds on $\alpha^{k,l}$ can be found in Proof of Theorem \ref{thm1}.
Conditions $(i)-(ii)$ and Lemma 4.7.2 in Puterman (1994) guarantee that $\alpha^{0,0}(\{s,s^L_2,s^L_1\})$ is isotone in $(\{s,s^L_2,s^L_1\})$. Hence,
\begin{align}
&\sum_{s^{L}_2 \neq z^{F_1} } \sigma_{s^{L}_2 ,z^{F_1}} \alpha^{0,0}(\{s,z^{F_1},s^L_1\})  \notag \\
&=\sum_{s^{L}_2 \neq z^{F_1} }\sigma_{s^{L}_2,z^{F_1}}R^L_\delta(\{s,z^{F_1},s^L_1\}) + \beta \sum_{s^{L}_2 \neq z^{F_1} }\sigma_{s^{L}_2,z^{F_1}} \sum_{s'}P(s'|s,\delta(\{s,z^{F_1},s^L_1\}))\alpha^{0,0}(\{s',s^{L'}_2,s^{L'}_1\}) \notag\\
&=\sum_{s^{L}_2 \neq z^{F_1} }\sigma_{s^{L}_2,z^{F_1}}R^L_\delta(\{s,z^{F_1},s^L_1\}) + \beta \sum_{s'}[\sum_{s^{L}_2 \neq z^{F_1} }\sigma_{s^{L}_2,z^{F_1}} P(s'|s,\delta(\{s,z^{F_1},s^L_1\}))]\alpha^{0,0}(\{s',s^{L'}_2,s^{L'}_1\}) \notag\\
&\leq \sum_{s^{L}_1 \neq z^{F_2} }\sigma_{s^{L}_1,z^{F_2}}R^L_\delta(\{s,s^L_1,z^{F_2}\}) + \beta \sum_{s'}[\sum_{s^{L}_1 \neq z^{F_2} }\sigma_{s^{L}_1,z^{F_2}} P(s'|s,\delta(\{s,s^L_2,z^{F_2}\}))]\alpha^{0,0}(\{s',s^{L'}_2,s^{L'}_1\}) \notag \\
&=\sum_{s^{L}_1 \neq z^{F_2} } \sigma_{s^{L}_1 ,z^{F_2}} \alpha^{0,0}(\{s,s^L_2,z^{F_2}\}) \notag
\end{align}
where inequality comes from Conditions (iii)-(iv), the isotonicity of $\alpha^{0,0}$ in $(\{s,s^L_2,s^L_1\})$, and Lemma 4.7.2 in Puterman (1994). Hence, $\Delta^{0,0} \leq \bar{\Delta}^{0,0}$. The result follows by the definition of $\alpha^{0,1}$ and $\alpha^{1,0}$ in terms of $\bar{\Delta}^{0,0}$ and $ \Delta^{0,0}$, respectively.

\textbf{Challenges in analyzing finite memory policies.} Note that by definition, both $\alpha^{1,0}$ and $\alpha^{0,1}$ can be represented in terms of $\alpha^{0,0}$. Hence, if $\alpha^{0,0}(\bar{\zeta}(s', z^{F_1'}, s^{L'}_1, \zeta, \delta(\zeta))) \leq \alpha^{0,0}(\bar{\zeta}(s',s^{L'}_2, z^{F_2'},  \zeta, \delta(\zeta))), \forall z^{F_1'}, z^{F_2'}$, then $\alpha^{1,0} \leq \alpha^{0,1}$. Conditions satisfying this inequality are indeed not hard to find. However, such conditions would be rather conservative. It would imply that the impact of having errors on follower 1 is always worse than that of having errors on follower 2, independent of the magnitude of the errors. It is more meaningful to compare the \textit{average} performances, $w=\sum_{s^{L'}_2 \neq z^{F'_1}} \sigma_{s^{L'}_2 ,z^{F_1'}}\alpha^{0,0}(\bar{\zeta}(s', z^{F_1'}, s^{L'}_1, \zeta, \delta(\zeta)))$ and $w'=\sum_{s^{L'}_1 \neq z^{F'_2}} \sigma_{s^{L'}_1 ,z^{F_2'}}\alpha^{0,0}(\bar{\zeta}(s', s^{L'}_2, z^{F_2'}, \zeta, \delta(\zeta)))$, as we did in Theorem \ref{quantify}. However, it is not straightforward for us to find a contraction operator $H$ satisfying $Hw=w$ without the zero-memory assumption.
\endproof

\section*{Appendix B: Parameters for the Examples}
\renewcommand{\labelitemi}{$\bullet$}
{\scriptsize
	Example 1 Meaning of each parameter:
	  \begin{table}[H]
		\centering
		\begin{tabular}{|c|c|l|}
			\hline
\multirow{5}{*}{State Spaces} & \multirow{3}{*}{Army $S^L = S^L_1\times S^L_2$} & $s^L=(s^L_1,s^L_2)$ where $s^L_i$ is the intensity \\ 
& & level of army activity in vicinity of\\
& &  position $i$\\
			\cline{2-3}
& \multirow{2}{*}{Militants $S^F$} & The state of militants is the allocation \\
&& of their heavy weapons. \\
\hline 	
\multirow{3}{*}{Observation Space} & \multirow{3}{*}{Militants $Z^{F_i}$} & Group $i$'s observation of intensity level \\
&&of army activity at location $j$, shared by \\
&&militant group $j$. Thus, $Z^{F_i}=S^L_j$. \\
\hline
	\multirow{3}{*}{Action Spaces} & Army $A^L$ & Patrol positions and kill militants \\ 
	\cline{2-3}
	& \multirow{2}{*}{Militants $A^F$} & Group $i$ will either use heavy weapons\\
	&& or reallocate their heavy weapons. \\
	\hline 	
	\multirow{8}{*}{Transition Probabilities} & \multirow{3}{*}{Army $P(s^{L'}|s,a)$} & The intensity level of army activity at the\\ 
	& & next decision epoch depends on the states \\
	& & and actions of all agents\\
	\cline{2-3}
	& \multirow{5}{*}{Militants $P(s^{F'}|s,a)$} & The allocation of limited heavy weapons \\
	&& next round depends on the current state\\
	&& and actions of all players. It also captures \\
	&& the agility of the militants (i.e., how fast\\
	&& they can move critical resources)\\
	\hline 		
\multirow{5}{*}{Reward Structures} & \multirow{2}{*}{Army $r^L(s,a)$} & Area reclaimed and the number of army  \\ 
& & troops killed by militants each round\\
\cline{2-3}
& \multirow{3}{*}{Militants $r^{F_i}(s,a)$} & The number of army troops killed each \\
&& round, impact of propaganda and symbolic\\
&&  purposes, etc. \\
\hline	
		\end{tabular}
	\end{table}
	Example 1 value of parameters - under a fixed leader policy:
	$$
	\begin{array}{cc}
	\mbox{Game 1: reward structure} & \mbox{Game 2: reward structure} \\
	r^L=(s,a^L=\mbox{fixed},a^F) \hspace{10pt}|\hspace{10pt} r^F=(s,a^L=\mbox{fixed},a^F) & r^L=(s,a^L=\mbox{fixed},a^F) \hspace{10pt}|\hspace{10pt} r^F=(s,a^L=\mbox{fixed},a^F)\\
	\left[ \begin{array}{cccc:ccccc}
	29&	28&	21&	30&	-29&	-28&	-21&	-30\\
	-24&	-33&	-13&	-32&	24&	33&	13&	32\\
	-41&	-30&	-8&	-24&	41&	30&	8&	24\\
	33&	17&	38&	46&	-33&	-17&	-38&	-46\\
	30&	-26&	46&	23&	-30&	26&	-46&	-23\\
	35&	22&	39&	28&	-35&	-22&	-39&	-28\\
	-46&	-8&	-5&	-22&	46&	8&	5&	22\\
	9&	-34&	38&	-3&	-9&	34&	-38&	3\\
	\end{array} \right] &
	\left[ \begin{array}{cccc:ccccc}
	29&	28&	21&	30&	-26&	-10&	-24&	-2\\
	-24&	-33&	-13&	-32&	4&	2&	24&	15\\
	-41&	-30&	-8&	-24&	18&	32&	48&	40\\
	33&	17&	38&	46&	-6&	-49&	-36&	-16\\
	30&	-26&	46&	23&	-16&	16&	-15&	-5\\
	35&	22&	39&	28&	-14&	-18&	-22&	-40\\
	-46&	-8&	-5&	-22&	40&	39&	13&	31\\
	9&	-34&	38&	-3&	-8&	49&	-44&	5\\
	\end{array} \right]
	\end{array}
	$$
	\begin{center}
		Transition Probabilities for both Game 1 and Game 2
	\end{center}
	$$	\begin{array}{cc}
	P(s'|s,a=(0^L, 0^{F_1}, 0^{F_2})) & P(s'|s,a=(0^L, 0^{F_1}, 1^{F_2})) \\
	\begin{bmatrix}
	0.17&	0.14&	0.02&	0.03&	0.16&	0.14&	0.18&	0.16\\
	0.25&	0.07&	0.08&	0.09&	0.16&	0.21&	0.14&	0\\
	0.26&	0.11&	0.01&	0.02&	0.25&	0.22&	0.07&	0.06\\
	0&	0.05&	0.04&	0.1&	0.06&	0.24&	0.22&	0.29\\
	0.22&	0.05&	0.23&	0.22&	0.18&	0.05&	0.02&	0.03\\
	0.16&	0.16&	0.13&	0&	0.09&	0.07&	0.2&	0.19\\
	0.12&	0.15&	0.1&	0.11&	0.22&	0.14&	0.09&	0.07\\
	0.1&	0.09&	0.18&	0.19&	0.17&	0.12&	0.05&	0.1
	\end{bmatrix} &
	\begin{bmatrix}
	0.33&	0.15&	0.14&	0.05&	0.22&	0.1&	0.01&	0\\
	0.22&	0.21&	0&	0.05&	0.06&	0.06&	0.16&	0.24\\
	0.18&	0.14&	0.14&	0.17&	0.04&	0.1&	0.19&	0.04\\
	0.19&	0.11&	0.08&	0.03&	0.22&	0.13&	0.12&	0.12\\
	0.18&	0.03&	0.13&	0.16&	0.11&	0.11&	0.11&	0.17\\
	0.03&	0.08&	0.14&	0.12&	0.16&	0.16&	0.13&	0.18\\
	0.01&	0.09&	0.15&	0.25&	0.12&	0.07&	0.13&	0.18\\
	0.21&	0.03&	0.02&	0.04&	0.27&	0.2&	0.05&	0.18\\
	\end{bmatrix}\\
	P(s'|s,a=(0^L, 1^{F_1}, 0^{F_2})) & P(s'|s,a=(0^L, 1^{F_1}, 1^{F_2})) \\
	\begin{bmatrix}
	0.18&	0.11&	0.16&	0&	0.11&	0.07&	0.16&	0.21\\
	0.16&	0.07&	0.28&	0.09&	0.1&	0.01&	0.09&	0.2\\
	0.14&	0.1&	0.15&	0.11&	0.17&	0.04&	0.13&	0.16\\
	0.03&	0.26&	0.23&	0.08&	0.13&	0.05&	0.11&	0.11\\
	0.15&	0.15&	0.08&	0.13&	0.16&	0.15&	0.1&	0.08\\
	0.08&	0.04&	0.05&	0.2&	0.06&	0.13&	0.22&	0.22\\
	0.11&	0.26&	0.09&	0.06&	0.22&	0.13&	0.05&	0.08\\
	0.11&	0.14&	0&	0.04&	0.36&	0.11&	0.04&	0.2\\
	\end{bmatrix} &
	\begin{bmatrix}
	0.09&	0.02&	0.1&	0.15&	0.2&	0.17&	0.03&	0.24\\
	0.03&	0.12&	0.24&	0.04&	0.06&	0.29&	0.04&	0.18\\
	0.07&	0.04&	0.28&	0.04&	0.24&	0.03&	0.05&	0.25\\
	0.03&	0.2&	0.16&	0.06&	0.31&	0.1&	0.02&	0.12\\
	0.2&	0.21&	0.09&	0.2&	0.01&	0.08&	0.15&	0.06\\
	0.27&	0.25&	0.05&	0.11&	0.03&	0.23&	0.02&	0.04\\
	0&	0.06&	0.19&	0.05&	0.22&	0.08&	0.19&	0.21\\
	0.16&	0.05&	0.03&	0.08&	0.19&	0.2&	0.17&	0.12\\
	\end{bmatrix}
	\end{array}$$
	
		Example 2 Meaning of each parameter:
	\begin{table}[H]
		\centering
		\begin{tabular}{|c|c|l|}
			\hline
			\multirow{5}{*}{State Spaces} & Army command  & $s^L=(s^L_1,s^L_2)$ where $s^L_i$ is the location \\ 
			&$S^L = S^L_1\times S^L_2$ & of combat unit $i$. \\
			\cline{2-3}
			& \multirow{3}{*}{Combat units $S^F$} & The state of the combat units is the  \\
			&& distribution of the heavy weapons of the   \\
			&& militants they are facing. \\
			\hline 	
			\multirow{2}{*}{Observation Space} & \multirow{2}{*}{Combat units $Z^{F_i}$} & Unit $i$'s observation of the location of \\
			&& unit $j$, shared by unit $j$. Thus, $Z^{F_i}=S^L_j$.\\
			\hline
			\multirow{4}{*}{Action Spaces} & Army command $A^L$ & Support combat units according to a plan. \\ 
			\cline{2-3}
			& \multirow{3}{*}{Combat units $A^F$} & Unit $i$ will either fight the militants \\
			&& opposing it or provide fire support\\
			&& to unit $j$. \\
			\hline 	
			\multirow{3}{*}{Transition Probability} & \multirow{3}{*}{ $P(s^{'}|s,a)$} & The state of next decision epoch depends\\ 
			& & on the current state and actions of \\
			& & all agents.\\
			\hline	
			\multirow{3}{*}{Reward Structures} & \multirow{2}{*}{Army command $r^L(s,a)$} & Area reclaimed and the number of army  \\ 
			& & troops killed by militants each round.\\
			\cline{2-3}
			& Combat units $r^{F_i}(s,a)$ & The number of militants killed per round. \\
			\hline	
		\end{tabular}
	\end{table}
	
	Example 2 parameters - under a fixed leader policy:
	$$
	\begin{array}{cc}
	\mbox{Game 1: reward structure} & \mbox{Game 2: reward structure} \\
	r^L=(s,a^L=\mbox{fixed},a^F) \hspace{10pt}|\hspace{10pt} r^F=(s,a^L=\mbox{fixed},a^F) & r^L=(s,a^L=\mbox{fixed},a^F) \hspace{10pt}|\hspace{10pt} r^F=(s,a^L=\mbox{fixed},a^F)\\
	\left[ \begin{array}{cccc:ccccc}
	21&	14&	-1&	-11&	21&	14&	-1&	-11\\
	-13&	47&	-20&	37&	-13&	47&	-20&	37\\
	-30&	-14&	-20&	-29&	-30&	-14&	-20&	-29\\
	31&	25&	-18&	-31&	31&	25&	-18&	-31\\
	-22&	-20&	-47&	-14&	-22&	-20&	-47&	-14\\
	48&	-42&	49&	28&	48&	-42&	49&	28\\
	49&	21&	10&	22&	49&	21&	10&	22\\
	-12&	-45&	-27&	-27&	-12&	-45&	-27&	-27\\
	\end{array} \right] &
	\left[ \begin{array}{cccc:ccccc}
	21&	14&	-1&	-11&	33&	49&	-25&	-48\\
	-13&	47&	-20&	37&	-29&	50&	-20&	3\\
	-30&	-14&	-20&	-29&	-38&	-13&	-50&	-48\\
	31&	25&	-18&	-31&	29&	32&	-38&	-46\\
	-22&	-20&	-47&	-14&	-45&	-24&	-2&	-38\\
	48&	-42&	49&	28&	23&	-33&	41&	45\\
	49&	21&	10&	22&	34&	25&	26&	31\\
	-12&	-45&	-27&	-27&	-39&	-28&	-14&	-46\\
	\end{array} \right]
	\end{array}
	$$
	\begin{center}
		Transition Probabilities for both Game 1 and Game 2
	\end{center}
	$$	\begin{array}{cc}
	P(s'|s,a=(0^L, 0^{F_1}, 0^{F_2})) & P(s'|s,a=(0^L, 0^{F_1}, 1^{F_2})) \\
	\begin{bmatrix}
	0&	0.06&	0.2&	0.18&	0.21&	0.18&	0.09&	0.08\\
	0.06&	0.18&	0.26&	0.02&	0.02&	0.13&	0.2&	0.13\\
	0.11&	0.2&	0.13&	0.01&	0.2&	0.04&	0.13&	0.18\\
	0.2&	0.14&	0.07&	0.02&	0.19&	0.1&	0.14&	0.14\\
	0.3&	0.13&	0.14&	0&	0.08&	0.21&	0.1&	0.04\\
	0.12&	0.08&	0.22&	0.21&	0.17&	0.06&	0.03&	0.11\\
	0.14&	0.05&	0.06&	0.1&	0.05&	0.21&	0.2&	0.19\\
	0.31&	0.17&	0.13&	0.07&	0.03&	0.22&	0.02&	0.05\\
	\end{bmatrix} &
	\begin{bmatrix}
	0.14&	0.17&	0.13&	0.01&	0.1&	0.07&	0.2&	0.18\\
	0.02&	0.29&	0.12&	0.01&	0.03&	0.16&	0.11&	0.26\\
	0.14&	0&	0.13&	0.18&	0.1&	0.18&	0.07&	0.2\\
	0.25&	0.03&	0.04&	0.05&	0.19&	0.18&	0.06&	0.2\\
	0.28&	0.1&	0.01&	0.08&	0.14&	0.06&	0.1&	0.23\\
	0.15&	0.07&	0.2&	0.05&	0.2&	0.19&	0.04&	0.1\\
	0.1&	0.1&	0.11&	0.04&	0.28&	0.13&	0.07&	0.17\\
	0.15&	0.15&	0.13&	0.12&	0.11&	0.13&	0.1&	0.11\\
	\end{bmatrix}\\
	P(s'|s,a=(0^L, 1^{F_1}, 0^{F_2})) & P(s'|s,a=(0^L, 1^{F_1}, 1^{F_2})) \\
	\begin{bmatrix}
	0.02&	0.14&	0.27&	0.18&	0.13&	0.03&	0&	0.23\\
	0.07&	0.06&	0.21&	0.15&	0.06&	0.18&	0.18&	0.09\\
	0.1&	0.17&	0.15&	0.03&	0.2&	0.03&	0.12&	0.2\\
	0.21&	0.09&	0.06&	0.15&	0.16&	0.02&	0.11&	0.2\\
	0.07&	0.11&	0.21&	0.03&	0.21&	0.12&	0.09&	0.16\\
	0.09&	0&	0.28&	0.27&	0.18&	0.07&	0.08&	0.03\\
	0.2&	0.18&	0.06&	0.2&	0.01&	0.04&	0.03&	0.28\\
	0.15&	0.12&	0.17&	0.07&	0.2&	0.1&	0.1&	0.09\\
	\end{bmatrix} &
	\begin{bmatrix}
	0.11&	0.12&	0.18&	0.17&	0.16&	0.08&	0.03&	0.15\\
	0.19&	0.03&	0.17&	0.15&	0&	0.27&	0.09&	0.1\\
	0.17&	0.03&	0.13&	0.2&	0.17&	0.05&	0.09&	0.16\\
	0.1&	0.05&	0.21&	0.1&	0.11&	0.21&	0.12&	0.1\\
	0.05&	0.15&	0.07&	0.17&	0.14&	0.16&	0.08&	0.18\\
	0.21&	0.11&	0.09&	0.06&	0.2&	0.07&	0.1&	0.16\\
	0.07&	0.15&	0.12&	0.16&	0.19&	0.15&	0.02&	0.14\\
	0.15&	0.16&	0.04&	0.05&	0.17&	0.09&	0.17&	0.17\\
	\end{bmatrix}
	\end{array}$$
	
	Example 3:
	\begin{center}
		reward structure $r^k(s^L_1,s^L_2,s^F,a^L=\mbox{fixed},a^{F_1},a^{F_2}), k \in \{L, F_1, F_2\}$
	\end{center}
	
	$$r^L=(s,a^L=\mbox{fixed},a^F) \hspace{10pt}|\hspace{10pt} r^F=(s,a^L=\mbox{fixed},a^F) $$
	$$\left[ \begin{array}{cccc:ccccc}
	17&	48&	48&	-5&	17&	48&	48&	-5\\
	-13&	-16&	-16&	-25&	-13&	-16&	-16&	-25\\
	-20&	-16&	-16&	-19&	-20&	-16&	-16&	-19\\
	-32&	5&	5&	37&	-32&	5&	5&	37\\
	-20&	-16&	-16&	-19&	-20&	-16&	-16&	-19\\
	-32&	5&	5&	37&	-32&	5&	5&	37\\
	-22&	40&	40&	27&	-22&	40&	40&	27\\
	11&	-28&	-28&	47&	11&	-28&	-28&	47\\
	\end{array} \right]$$
	\begin{center}
		Transition Probabilities
	\end{center}
	$$	\begin{array}{cc}
	P(s'|s,a=(0^L, 0^{F_1}, 0^{F_2})) & P(s'|s,a=(0^L, 0^{F_1}, 1^{F_2})) \\
	\begin{bmatrix}
	0.06&	0.05&	0.12&	0.14&	0.12&	0.14&	0.17&	0.2\\
	0.17&	0.16&	0.22&	0.07&	0.21&	0.07&	0.05&	0.05\\
	0.31&	0.04&	0.01&	0.22&	0.01&	0.22&	0.04&	0.15\\
	0.13&	0.15&	0.07&	0.2&	0.07&	0.2&	0.07&	0.11\\
	0.31&	0.04&	0.01&	0.22&	0.01&	0.22&	0.04&	0.15\\
	0.13&	0.15&	0.07&	0.2&	0.07&	0.2&	0.07&	0.11\\
	0.33&	0.04&	0.15&	0.02&	0.15&	0.02&	0.27&	0.02\\
	0.19&	0.15&	0.12&	0.08&	0.12&	0.08&	0.21&	0.05\\
	\end{bmatrix} &
	\begin{bmatrix}
	0.2&	0.18&	0.03&	0.01&	0.03&	0.01&	0.24&	0.3\\
	0.25&	0.1&	0.04&	0.15&	0.04&	0.15&	0.11&	0.16\\
	0.07&	0.13&	0.14&	0.16&	0.14&	0.15&	0.13&	0.08\\
	0.31&	0.12&	0.07&	0.03&	0.07&	0.03&	0.09&	0.28\\
	0.07&	0.13&	0.14&	0.15&	0.14&	0.16&	0.13&	0.08\\
	0.31&	0.12&	0.07&	0.03&	0.07&	0.03&	0.09&	0.28\\
	0.24&	0.18&	0.01&	0.17&	0.01&	0.17&	0.18&	0.04\\
	0.18&	0.02&	0.1&	0.19&	0.1&	0.19&	0.18&	0.04\\
	\end{bmatrix}\\
	P(s'|s,a=(0^L, 1^{F_1}, 0^{F_2})) & P(s'|s,a=(0^L, 1^{F_1}, 1^{F_2})) \\
	\begin{bmatrix}
	0.2&	0.18&	0.03&	0.01&	0.03&	0.01&	0.24&	0.3\\
	0.25&	0.1&	0.04&	0.15&	0.04&	0.15&	0.11&	0.16\\
	0.07&	0.13&	0.14&	0.16&	0.14&	0.15&	0.13&	0.08\\
	0.31&	0.12&	0.07&	0.03&	0.07&	0.03&	0.09&	0.28\\
	0.07&	0.13&	0.14&	0.15&	0.14&	0.16&	0.13&	0.08\\
	0.31&	0.12&	0.07&	0.03&	0.07&	0.03&	0.09&	0.28\\
	0.24&	0.18&	0.01&	0.17&	0.01&	0.17&	0.18&	0.04\\
	0.18&	0.02&	0.1&	0.19&	0.1&	0.19&	0.18&	0.04\\
	\end{bmatrix} &
	\begin{bmatrix}
	0.06&	0.05&	0.1&	0.1&	0.1&	0.1&	0.18&	0.31\\
	0.01&	0.38&	0.14&	0.07&	0.14&	0.07&	0.02&	0.17\\
	0.19&	0.24&	0.11&	0.06&	0.11&	0.06&	0.05&	0.18\\
	0.02&	0.25&	0.22&	0.01&	0.22&	0.01&	0.06&	0.21\\
	0.19&	0.24&	0.11&	0.06&	0.11&	0.06&	0.05&	0.18\\
	0.02&	0.25&	0.22&	0.01&	0.22&	0.01&	0.06&	0.21\\
	0.1&	0.02&	0.18&	0.16&	0.18&	0.16&	0.07&	0.13\\
	0.04&	0.14&	0.14&	0.13&	0.14&	0.13&	0.14&	0.14\\
	\end{bmatrix}
	\end{array}$$
	Example 4:
	\begin{center}
		reward structure $r^k(s^L_1,s^L_2,s^F,a^L=\mbox{fixed},a^{F_1},a^{F_2}), k \in \{L, F_1, F_2\}$
	\end{center}
	$$r^L=(s,a^L=\mbox{fixed},a^F) \hspace{10pt}|\hspace{10pt} r^F=(s,a^L=\mbox{fixed},a^F) $$
	$$\left[ \begin{array}{cccc:ccccc}
	-1028&	-968&	-968&	-968&	-1028&	-968&	-968&	-968\\
	-962&	-1006&	-1006&	-952&	-962&	-1006&	-1006&	-952\\
	-964&	-1038&	-1038&	-1029&	-964&	-1038&	-1038&	-1029\\
	-958&	-967&	-967&	-987&	-958&	-967&	-967&	-987\\
	-964&	-1038&	-1038&	-1029&	-964&	-1038&	-1038&	-1029\\
	-958&	-967&	-967&	-987&	-958&	-967&	-967&	-987\\
	1025&	1004&	1004&	998&	1025&	1004&	1004&	998\\
	965&	1009&	1009&	1036&	965&	1009&	1009&	1036\\
	\end{array} \right]$$
	\begin{center}
		Transition Probabilities
	\end{center}
	$$	\begin{array}{cc}
	P(s'|s,a=(0^L, 0^{F_1}, 0^{F_2})) & P(s'|s,a=(0^L, 0^{F_1}, 1^{F_2})) \\
	\begin{bmatrix}
	0.11&	0.03&	0&	0.17&	0&	0.17&	0.18&	0.34\\
	0.08&	0.2&	0.05&	0&	0.05&	0&	0.08&	0.54\\
	0.12&	0.15&	0.12&	0&	0.12&	0&	0.36&	0.13\\
	0.18&	0.05&	0.02&	0.08&	0.02&	0.08&	0.17&	0.4\\
	0.12&	0.15&	0.12&	0&	0.12&	0&	0.36&	0.13\\
	0.18&	0.05&	0.02&	0.08&	0.02&	0.08&	0.17&	0.4\\
	0.16&	0.09&	0.03&	0&	0.03&	0&	0.59&	0.1\\
	0.1&	0.12&	0.02&	0.1&	0.02&	0.1&	0.42&	0.12\\
	\end{bmatrix} &
	\begin{bmatrix}
	0.09&	0.1&	0.01&	0.09&	0.01&	0.09&	0.35&	0.26\\
	0.11&	0.15&	0&	0.05&	0&	0.05&	0.15&	0.49\\
	0.07&	0.1&	0.01&	0.1&	0.01&	0.1&	0.14&	0.47\\
	0.03&	0.14&	0.08&	0.03&	0.08&	0.03&	0.34&	0.27\\
	0.07&	0.1&	0.01&	0.1&	0.01&	0.1&	0.14&	0.47\\
	0.03&	0.14&	0.08&	0.03&	0.08&	0.03&	0.34&	0.27\\
	0.05&	0.16&	0.04&	0.11&	0.04&	0.11&	0.25&	0.24\\
	0.12&	0.06&	0.04&	0.09&	0.04&	0.09&	0.46&	0.1\\
	\end{bmatrix}\\
	P(s'|s,a=(0^L, 1^{F_1}, 0^{F_2})) & P(s'|s,a=(0^L, 1^{F_1}, 1^{F_2})) \\
	\begin{bmatrix}
	0.09&	0.1&	0.01&	0.09&	0.01&	0.09&	0.35&	0.26\\
	0.11&	0.15&	0&	0.05&	0&	0.05&	0.15&	0.49\\
	0.07&	0.1&	0.01&	0.1&	0.01&	0.1&	0.14&	0.47\\
	0.03&	0.14&	0.08&	0.03&	0.08&	0.03&	0.34&	0.27\\
	0.07&	0.1&	0.01&	0.1&	0.01&	0.1&	0.14&	0.47\\
	0.03&	0.14&	0.08&	0.03&	0.08&	0.03&	0.34&	0.27\\
	0.05&	0.16&	0.04&	0.11&	0.04&	0.11&	0.25&	0.24\\
	0.12&	0.06&	0.04&	0.09&	0.04&	0.09&	0.46&	0.1\\
	\end{bmatrix} &
	\begin{bmatrix}
	0.06&	0.09&	0.06&	0.09&	0.06&	0.09&	0.54&	0.01\\
	0.08&	0.05&	0.09&	0.04&	0.09&	0.04&	0.35&	0.26\\
	0.06&	0.12&	0.06&	0.06&	0.06&	0.06&	0.23&	0.35\\
	0.06&	0.04&	0.06&	0.09&	0.06&	0.09&	0.37&	0.23\\
	0.06&	0.12&	0.06&	0.06&	0.06&	0.06&	0.23&	0.35\\
	0.06&	0.04&	0.06&	0.09&	0.06&	0.09&	0.37&	0.23\\
	0.02&	0.03&	0.06&	0.09&	0.06&	0.09&	0.35&	0.3\\
	0&	0.02&	0.07&	0.12&	0.07&	0.12&	0.1&	0.5\\
	\end{bmatrix}
	\end{array}$$
	Example in table \ref{alphacompare}:
	\begin{center}
		Reward Structure
	\end{center}
	
	\begin{table}[H]
		\begin{center}
			{\scriptsize
				\begin{tabular}{ | c | c | c | c | c | c | c | c |} \hline			

					$\zeta$ &$r_\delta(\zeta)$ & $\zeta$ &$r_\delta(\zeta)$ &$\zeta$ &$r_\delta(\zeta)$ &$\zeta$ &$r_\delta(\zeta)$  \\
					\hline
					[0 0] [0 0] 0	&	10	&	[0 0] [0 1] 0	&	26	&	[1 0] [0 0] 0	&	48	&	[1 0] [0 1] 0	&	65	\\\hline
					[0 1] [0 0] 0	&	6	&	[0 1] [0 1] 0	&	30	&	[1 1] [0 0] 0	&	45	&	[1 1] [0 1] 0	&	68	\\\hline
					[0 0] [1 0] 0	&	8	&	[0 0] [1 1] 0	&	25	&	[1 0] [1 0] 0	&	50	&	[1 0] [1 1] 0	&	66	\\\hline
					[0 1] [1 0] 0	&	5	&	[0 1] [1 1] 0	&	28	&	[1 1] [1 0] 0	&	46	&	[1 1] [1 1] 0	&	70	\\\hline
					[0 0] [0 0] 1	&	20	&	[0 0] [0 1] 1	&	36	&	[1 0] [0 0] 1	&	58	&	[1 0] [0 1] 1	&	75	\\\hline
					[0 1] [0 0] 1	&	16	&	[0 1] [0 1] 1	&	40	&	[1 1] [0 0] 1	&	55	&	[1 1] [0 1] 1	&	78	\\\hline
					[0 0] [1 0] 1	&	18	&	[0 0] [1 1] 1	&	35	&	[1 0] [1 0] 1	&	60	&	[1 0] [1 1] 1	&	76	\\\hline
					[0 1] [1 0] 1	&	5	&	[0 1] [1 1] 1	&	38	&	[1 1] [1 0] 1	&	56	&	[1 1] [1 1] 1	&	80			
					\\\hline				
			\end{tabular}}	
		\end{center}
	\end{table}
	\begin{center}
		Transition Probabilities $P_\delta(\zeta) = P(s'|s,\delta(\zeta))$
	\end{center}
	$$\left[ \begin{array}{cccccccc:ccccccccc}
	0.7	&	0.05	&	0.05	&	0.05	&	0.05	&	0.05	&	0.05	&	0	&	0.3	&	0.15	&	0.05	&	0.05	&	0.05	&	0.05	&	0.15	&	0.2	\\
	0.76	&	0.02	&	0.05	&	0.05	&	0.05	&	0.05	&	0.02	&	0	&	0.35	&	0.1	&	0.03	&	0.02	&	0.05	&	0.15	&	0.1	&	0.2	\\
	0.74	&	0.03	&	0.05	&	0.05	&	0.05	&	0.05	&	0.03	&	0	&	0.5	&	0.05	&	0.05	&	0.05	&	0.05	&	0.05	&	0.05	&	0.2	\\
	0.1	&	0.1	&	0	&	0.1	&	0.3	&	0.1	&	0.1	&	0.2	&	0.24	&	0.18	&	0.05	&	0.05	&	0.05	&	0.05	&	0.18	&	0.2	\\
	0.65	&	0.05	&	0.05	&	0.05	&	0.05	&	0.05	&	0.05	&	0.05	&	0.15	&	0.2	&	0.05	&	0.05	&	0.05	&	0.05	&	0.2	&	0.25	\\
	0.67	&	0.04	&	0.05	&	0.05	&	0.05	&	0.05	&	0.04	&	0.05	&	0.39	&	0.1	&	0.01	&	0.03	&	0.02	&	0.15	&	0.1	&	0.2	\\
	0.63	&	0.06	&	0.05	&	0.05	&	0.05	&	0.05	&	0.06	&	0.05	&	0.45	&	0.05	&	0.05	&	0.05	&	0.05	&	0.05	&	0.05	&	0.25	\\
	0	&	0.1	&	0.2	&	0.3	&	0.1	&	0	&	0.1	&	0.2	&	0.21	&	0.17	&	0.05	&	0.05	&	0.05	&	0.05	&	0.17	&	0.25	\\
	0.58	&	0.06	&	0.05	&	0.05	&	0.05	&	0.05	&	0.06	&	0.1	&	0.42	&	0.1	&	0.05	&	0.02	&	0.01	&	0.1	&	0.1	&	0.2	\\
	0.6	&	0.05	&	0.05	&	0.05	&	0.05	&	0.05	&	0.05	&	0.1	&	0.06	&	0.22	&	0.05	&	0.05	&	0.05	&	0.05	&	0.22	&	0.3	\\
	0.1	&	0.1	&	0	&	0	&	0.3	&	0.2	&	0.1	&	0.2	&	0.2	&	0.15	&	0.05	&	0.05	&	0.05	&	0.05	&	0.15	&	0.3	\\
	0.5	&	0.1	&	0.05	&	0.05	&	0.05	&	0.05	&	0.1	&	0.1	&	0.4	&	0.05	&	0.05	&	0.05	&	0.05	&	0.05	&	0.05	&	0.3	\\
	0.49	&	0.08	&	0.05	&	0.05	&	0.05	&	0.05	&	0.08	&	0.15	&	0.3	&	0.1	&	0.2	&	0.03	&	0.05	&	0.02	&	0.1	&	0.2	\\
	0.55	&	0.05	&	0.05	&	0.05	&	0.05	&	0.05	&	0.05	&	0.15	&	0.15	&	0.05	&	0.05	&	0.05	&	0.05	&	0.05	&	0.25	&	0.35	\\
	0.13	&	0.1	&	0.02	&	0.05	&	0	&	0.4	&	0.1	&	0.2	&	0.18	&	0.05	&	0.05	&	0.05	&	0.05	&	0.05	&	0.22	&	0.35	\\
	0.41	&	0.12	&	0.05	&	0.05	&	0.05	&	0.05	&	0.12	&	0.15	&	0.35	&	0.05	&	0.05	&	0.05	&	0.05	&	0.05	&	0.05	&	0.35	\\	
	\end{array} \right]$$
	%
}

\end{document}